\pdfoutput=1
\documentclass[11pt,a4paper]{amsart}

\usepackage{amsfonts, amsmath, amssymb, amsthm}
\usepackage{txfonts} %
\usepackage{latexsym}
\usepackage{graphicx}
\usepackage{url}
\usepackage{overpic}
\usepackage{color}
\usepackage{xspace}

\theoremstyle{plain}
\newtheorem{thm}{Theorem}

\newtheorem{prop}[thm]{Proposition}

\theoremstyle{definition} \newtheorem{defn}[thm]{Definition}
\newtheorem{ex}[thm]{Example}

\newtheorem{rem}[thm]{Remark}
\providecommand\ZZ{{\mathbb Z}}
\providecommand\RR{{\mathbb R}}

\providecommand\cF{{\mathcal F}}
\providecommand\cT{{\mathcal T}}
\providecommand\schlegel[2]{{\mathcal D}\left(#1,#2\right)}
\DeclareMathOperator{\conv}{conv}

\DeclareMathOperator{\Sym}{Sym}
\providecommand\SetOf[2]{\left\{#1\vphantom{#2}\,\right.\left|\,\vphantom{#1}#2\right\}}
\providecommand\smallSetOf[2]{\{#1\,|\,#2\}}

\providecommand\polymake{\texttt{polymake}\xspace}
\providecommand\JavaView{\texttt{JavaView}\xspace}
\providecommand\JReality{\texttt{JReality}\xspace}

\begin{document} \title[Drawing polytopal graphs]{Drawing polytopal graphs with \polymake}

\author[Gawrilow, Joswig, R\"orig, and Witte] {Ewgenij Gawrilow \and Michael Joswig \and
Thilo R\"orig \and Nikolaus Witte}
\address{Ewgenij Gawrilow, Institut f\"ur Mathematik,
MA 6-1, TU Berlin, 10623 Berlin, Germany}
\email{gawrilow@math.tu-berlin.de}
\address{Michael Joswig, Fachbereich Mathematik, AG~7, TU Darmstadt, 64289 Darmstadt,
Germany}
\email{joswig@mathematik.tu-darmstadt.de}
\address{Thilo R\"orig, Institut f\"ur
Mathematik, MA 6-2, TU Berlin, 10623 Berlin, Germany}
\email{thilosch@math.tu-berlin.de}
\address{Nikolaus Witte, Institut f\"ur Mathematik, MA 6-2, TU Berlin, 10623 Berlin,
Germany}
\email{witte@math.tu-berlin.de}

\date{November 15, 2007}

\begin{abstract} This note wants to explain how to obtain meaningful pictures of (possibly
high-dimensional) convex polytopes, triangulated manifolds, and other objects from the
realm of geometric combinatorics such as tight spans of finite metric spaces and tropical
polytopes.  In all our cases we arrive at specific, geometrically motivated, graph drawing
problems.  The methods displayed are implemented in the software system \polymake.
\end{abstract}

\maketitle

\section{Introduction}

\noindent Clearly, visualization is a key tool for doing experimental mathematics.
However, what to do if the geometric objects that we want to understand do not admit a
straightforward, natural way of visualization?  Reasons for this may include some of the
following: The objects live in Euclidean space but only of high dimension. The objects do
not come with any embedding into Euclidean space (or any other known geometry). We are
only interested in combinatorial features and thus are after a visualization which
abstracts from geometric ``randomness.''

In order to be able to give answers to the question posed we restrict our attention to
convex polytopes, triangulated surfaces and some other objects derived from them.  The
common theme will be that we will assign a graph to the object in question which then asks
for a suitable, meaningful visualization.  In some cases there is an obvious candidate for
such a graph like the vertex-edge graph of a polytope; see
Section~\ref{sec:polytope_graphs}.  In other cases there are interesting choices to be
made as for abstract simplicial manifolds; see Section~\ref{sec:pd_graphs}.

A graph~$G=(V,E)$ is a pair of a \emph{node set} $V$ and an \emph{edge set} $E$, consisting
of $2$-element subsets of $V$.  That is, our graphs are usually undirected and simple
(loop-less and without multiple edges).  We use the term \emph{polytopal graph} in a loose
sense: a polytopal graph is associated with some polytope(s) in one way or another.

It is a basic fact that each graph can be drawn in $\RR^3$ with straight edges and without
self-intersections; see Remark~\ref{rem:cyclic}.  An obvious question is how to find a
``good'' drawing representing a given graph (which may even admit self-intersections).
But these quality parameters depend on the context, and so we will discuss them, time and
again, in the subsequent sections.  The methods which we present are implemented in the
software package \polymake~\cite{polymake,gawrilow_joswig:POLYMAKE_I}.

The following is to describe the contents.  We start out with a very brief introduction to
convex polytopes (the objects \polymake primarily is designed for) and their graphs in
Section~\ref{sec:polytope_graphs}.  The situation for polytopes in dimension $\le 4$ is
special in that there is a canonical way of visualizing.  For dimensions $\le 3$ this is
obvious, and in dimension $4$ \emph{Schlegel diagrams} come in handy.  Their interactive
construction is the topic of Section~\ref{sec:schlegel_diagrams}.
Section~\ref{sec:pp_models} discusses the use of pseudo-physical models to produce
drawings of a graph.  A dynamic process is modeled which often converges to an acceptable
drawing of the graph.  Additionally, we discuss and exemplify the use of (pseudo-physical)
forces which turn out to be particularly useful for polytopes. In
Section~\ref{sec:rubber_bands} the rubber band method by
Maxwell~\cite{maxwell1864:on_reciprocal_figures} and Cremona~\cite{Cremona1890} is
explained, which expresses the dynamic process from a special planar pseudo-physical model
in terms of linear algebra.

We conclude this paper by discussing visualizations of polytopal graphs which come from
areas which recently received some attention.  This includes tight spans of finite metric
spaces (used for phylogenetic reconstructions in computational biology)
\cite{MR858908,MR2097310} and tropical polytopes (arising in combinatorial aspects of
algebraic geometry) \cite{DevelinSturmfels04,Joswig05,BlockYu06}.  Finally we introduce
\emph{pd-graphs} to visualize simplicial manifolds.

For general references to graph drawing see~\cite{TollisEtAl99,JuengerMutzel04}.

\section{Facts About Convex Polytopes and Their Graphs}
\label{sec:polytope_graphs}

\noindent A \emph{(convex) polyhedron} is the intersection of finitely many closed affine
half-spaces in some Euclidean space.  It is called a \emph{polytope} if it is bounded.
Each polytope is the convex hull of finitely many points, and vice versa.  Likewise, each
polyhedron can be described as the Minkowski sum of polytope, a finitely generated pointed
cone, and an affine subspace.  The polyhedron is \emph{pointed} if it does not contain any
affine subspace.  In the following we will assume that the polyhedron $P\subset\RR^d$ will
affinely span the space.  This is not much of a restriction since otherwise we can
continue our discussion by considering the affine span of~$P$ as the surrounding space.

A \emph{proper face} of such a full-dimensional polyhedron $P$ is the intersection with a
supporting affine hyperplane.  The empty set and $P$ itself are the two \emph{non-proper
faces}.  Faces are again polyhedra, and the set $\cF(P)$ of all faces is partially ordered
by inclusion.  This poset is naturally ranked by the \emph{dimension} of a face, which is
the dimension of its affine span.  The faces of dimensions $0$, $1$, $d-2$, $d-1$ are
called \emph{vertices}, \emph{edges}, \emph{ridges}, and \emph{facets}, respectively.

In the following let us assume that the polyhedron $P$ is bounded, that is, it is a
polytope.  If $P$ is an arbitrary pointed polyhedron, this again is not much of a
restriction: In this case $P$ is the image of a polytope under a projective linear
transformation.  If $P$ is not pointed then its projection to the orthogonal complement to
the \emph{lineality space} (that is, the maximal affine subspace that it contains) is
pointed.  Thanks to cone polarity the polytope $P$ has a polar polytope $P^*$ (with
respect to a chosen interior point).  Its face poset is the same as $\cF(P)$ but with
reversed inclusion.  In fact, if we take `intersection' as the meet-operation and `joint
convex hull' as the join-operation, the poset~$\cF(P)$ becomes an atomic and co-atomic
Eulerian lattice.

The \emph{(vertex-edge) graph} of $P$, denoted as $\Gamma(P)$, has the vertices of $P$ as
its nodes and the edges of $P$, well, as its edges.  There is also the dual graph
$\Gamma^*(P)$ formed from the facets and ridges.  We have $\Gamma^*(P)\cong\Gamma(P^*)$.
Notice that the Simplex Method from linear optimization walks along the vertex-edge graph
of a polyhedron (that is, the set of admissible points) to a vertex which is optimal with
respect to a given linear objective function.

A $0$-dimensional polytope is a point, a $1$-dimensional polytope is an edge, and a
$2$-dimensional polytope is a convex $n$-gon.  Hence, in terms of combinatorial properties
of their graphs or their face lattices, polytopes become interesting in dimension~$3$ and
beyond.  We list some facts known about polytopal graphs.

A graph is \emph{connected} if any two of its nodes are joined by an edge path.  We
call a graph with at least $k+1$ nodes \emph{$k$-connected} if removing any set of at most
$k-1$ nodes and the incident edges leaves the graph induced on the remaining nodes
connected.  This means that a $1$-connected graph is the same as a connected graph.  Note
that this must not be confused with the notion of higher connectivity common in topology.

\begin{thm}[{Steinitz~\cite{Steinitz22,SteinitzRademacher76}}]\label{thm:Steinitz} A graph
is isomorphic to the graph of a $3$-di\-men\-si\-o\-nal polytope if and only if it is
planar and $3$-connected.
\end{thm}

Here a graph is called \emph{planar} if it admits a drawing in $\RR^2$ without
self-in\-ter\-sec\-tions.

\begin{thm}[{Balinski~\cite{Balinski61}}] \label{thm:Balinski} The graph of a
$d$-dimensional polytope is \emph{$d$-connected}, in particular, the degree of each vertex
is at least~$d$.
\end{thm}

While graphs of polytopes do have a variety of properties which make them special,
altogether the class of polytopal graphs is large in the following sense: Kaibel and
Schwartz proved that the graph isomorphism problem restricted to polytopal graphs is as
hard as the general problem~\cite{KaibelSchwartz03}.

There is a particularly interesting class of polytopes: A $d$-dimensional polytope is
\emph{simple} if its graph is \emph{$d$-regular}, that is, each vertex has degree~$d$.  A
polytope is simple if and only if its polar is \emph{simplicial}, that is, each proper
face of the polar is a simplex.  In particular, if the affine hyperplanes spanned by the
facets of~$P$ are in general position, then $P$ is simple.  Examples of simple polytopes
include the $n$-gons, the simplices and the cubes (of arbitrary dimension).

\begin{thm}[{Blind and Mani~\cite{BlindMani87}}] \label{thm:BlindMani} Two simple
polytopes with isomorphic graphs have isomorphic face lattices.
\end{thm}

We will read this last result in the following way: For the visualization of
high-dimensional \emph{simple} polytopes it may suffice to visualize their vertex-edge
graphs.

For further information about polytopes the reader is referred to
Ziegler~\cite{ziegler:LOP}.

\section{Schlegel Diagrams}\label{sec:schlegel_diagrams}

\noindent The combinatorial features of a polytope are properties of its boundary, which
is of one dimension lower.  This can be exploited for visualization.  The Schlegel diagram
is a particular projection of the polytope onto one of its facets which still contains all
the combinatorial data.  This works in any dimension, but it is particularly useful to
understand the structure of $4$-dimensional polytopes by their $3$-dimensional
projections.

\begin{defn} Let $P \subset \RR^d$ be a $d$-polytope and $F$ one of its facets. Let $H_F =
\smallSetOf{ x \in \RR^d }{ ax = b }$ with $a \in \RR^d$ and $b \in \RR$ be the supporting
hyperplane of $F$, that is, $F = P \cap H_F$ and~$P$ is contained in the positive
half-space~$\smallSetOf{ x \in \RR^d }{ ax \leq b }$ defined by~$F$. Choose a point~$v$
contained in the negative half-space of~$F$ and in the positive half-spaces of all other
supporting hyperplanes of~$P$. We call~$v$ a point \emph{beyond}~$F$.  For every point $x
\in P$ we define the projection of~$x$ onto~$F$ by
  \[ \pi(x) \ = \ v + \frac{b - av}{ax-av}(x-v) \, .
  \] The \emph{Schlegel diagram of $P$ on $F$} is the polytopal subdivision of $F$
consisting of the image under $\pi$ of all proper faces of $P$ except $F$, that is,
  \[ \schlegel{P}{F} \ = \ \SetOf{ \pi(G) \subset F }{ G \in \cF(P)\setminus \{ \emptyset,
P,F \} } \, .
  \]
\end{defn}

The construction guarantees that a Schlegel diagram $\schlegel{P}{F}$ properly shows the
relative positions of the facets of $P$ and their intersections. For Schlegel diagrams of
$4$-polytopes it is common to visualize only their $1$-skeleta, that is, their graphs, in
order to be able to look ``inside'' $F$.  Schlegel diagrams of $2$- and $3$-dimensional
polytopes are shown in the subsequent Section \ref{subsec:all_schlegel}.

\begin{figure}\centering
  \includegraphics[width=.46\textwidth]{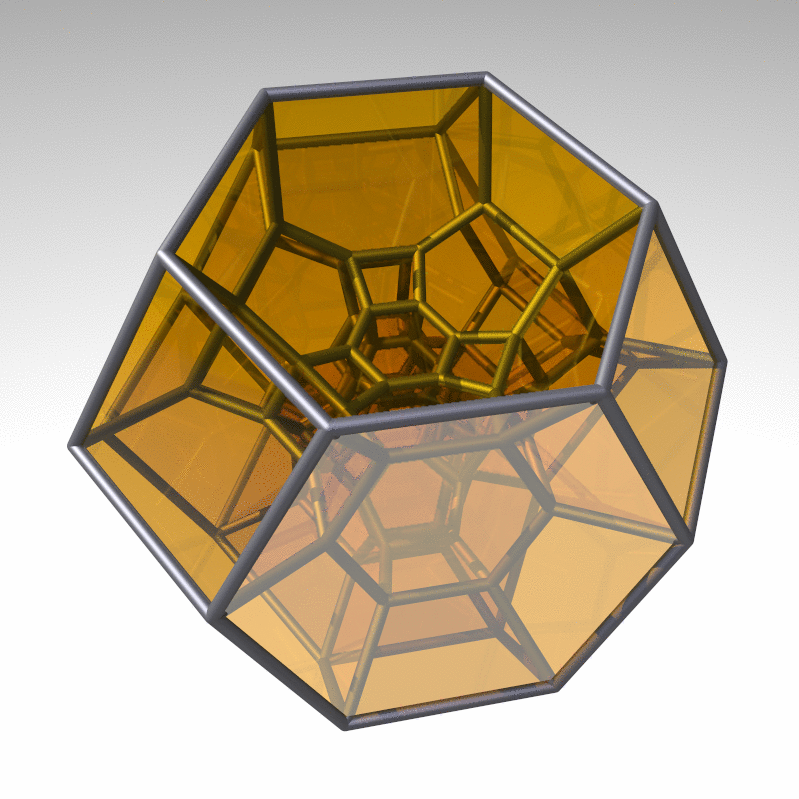} \hspace{.6cm}
  \includegraphics[width=.46\textwidth]{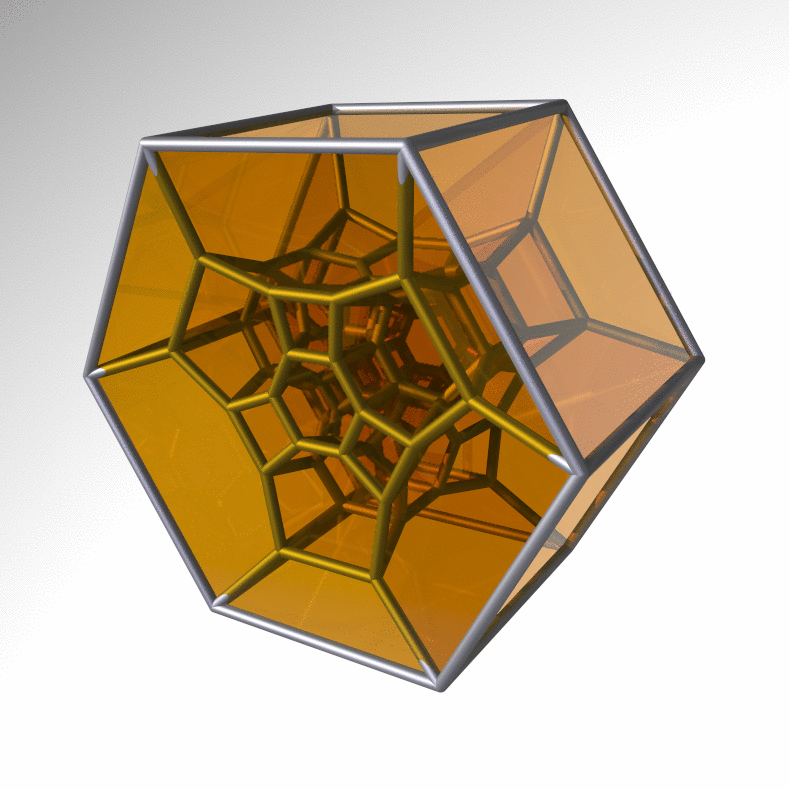}
  \caption{Schlegel diagrams of the $4$-permutohedron with respect to two different
facets: a $3$-permutohedron (left) and a prism over a hexagon (right).}
  \label{fig:schlegel:first}
\end{figure}

\begin{ex}\label{ex:permutohedron} If we write the permutation
  $\sigma\in\Sym\{1,2,\dots,n\}$ as the vector
  $v_\sigma=(\sigma(1),\sigma(2),\dots,\sigma(n))\in\RR^n$, then the \emph{permutohedron}
  of degree $n$ is the polytope
  \[ \Pi_{n-1} \ = \ \conv \SetOf{v_\sigma}{\sigma\in\Sym\{1,2,\dots,n\}} \, .
  \] The $(n-1)$-permutohedron $\Pi_{n-1}$ is a simple polytope of dimension $n-1$.  All
  its faces are products of lower-dimensional permutohedra.  The $2$-permutohedron is a
  hexagon, the $3$-permutohedron is an Archimedean solid with $14$ facets, all of which
  are squares or hexagons.  The facets of the $4$-permutohedron are $3$-permutohedra or
  prisms over hexagons.  Its two types of Schlegel diagrams are shown in
  Figure~\ref{fig:schlegel:first}.
\end{ex}

\begin{rem} The Schlegel diagram construction shows that the vertex-edge graphs of
$3$-polytopes are planar; see Theorem~\ref{thm:Steinitz}.
\end{rem}

\begin{rem}\label{rem:cyclic} Consider the moment curve $m:t\mapsto (t,t^2,t^3,t^4)$ in
$\RR^4$.  Taking $n$ distinct values $t_1<t_2<\dots<t_n$, the convex hull
$\conv\{m(t_1),m(t_2),\dots,m(t_n)\}$ of the corresponding points is a \emph{cyclic}
$4$-polytope on $n$ vertices.  As a special feature the vertex-edge graph of a cyclic
$4$-polytope is isomorphic to the complete graph.  Since each finite graph is a subgraph
of a complete graph, the Schlegel diagram construction establishes that each finite graph
admits a drawing in $\RR^3$ without self-intersections.  The cyclic polytopes are
simplicial.
\end{rem}

Since the geometry of the projection facet is preserved under the projection we cannot
expect to get a ``nice'' Schlegel diagram right away.  For some polytopes there may be a
good choice of a facet and a point beyond, but this choice may not be obvious.  The next
section discusses how to find a good Schlegel diagram interactively.

\subsection{Obtaining all Schlegel Diagrams of a Polytope}\label{subsec:all_schlegel}

\noindent The Schlegel diagram of a polytope $P$ depends on the projection facet $F$ and
the point~$v$ beyond~$F$, the \emph{viewpoint}. To obtain \emph{all} Schlegel diagrams of
a $4$-polytope we start with \emph{any} Schlegel diagram and describe an interactive
method to choose a different projection facet and a different point beyond.  This is
implemented in \polymake's interface to \JavaView \cite{javaview} and \JReality \cite{JReality}.

Since each facet is uniquely determined by its vertex set and all vertices are visible in
the Schlegel diagram, we are able to select the vertex set of another facet in the
Schlegel diagram $\schlegel{P}{F}$.  In fact, it suffices to mark sufficiently many
vertices, such that a unique facet containing them remains.  In the \JavaView and
\JReality graphical interfaces the user can mark points, and then press a button to get a
new window with a Schlegel diagram with respect to the facet defined by the marked points.
An error is issued if the facet is not uniquely specified.

\begin{figure}[thbp]
  \centering
  \begin{overpic}[width=6.5cm]{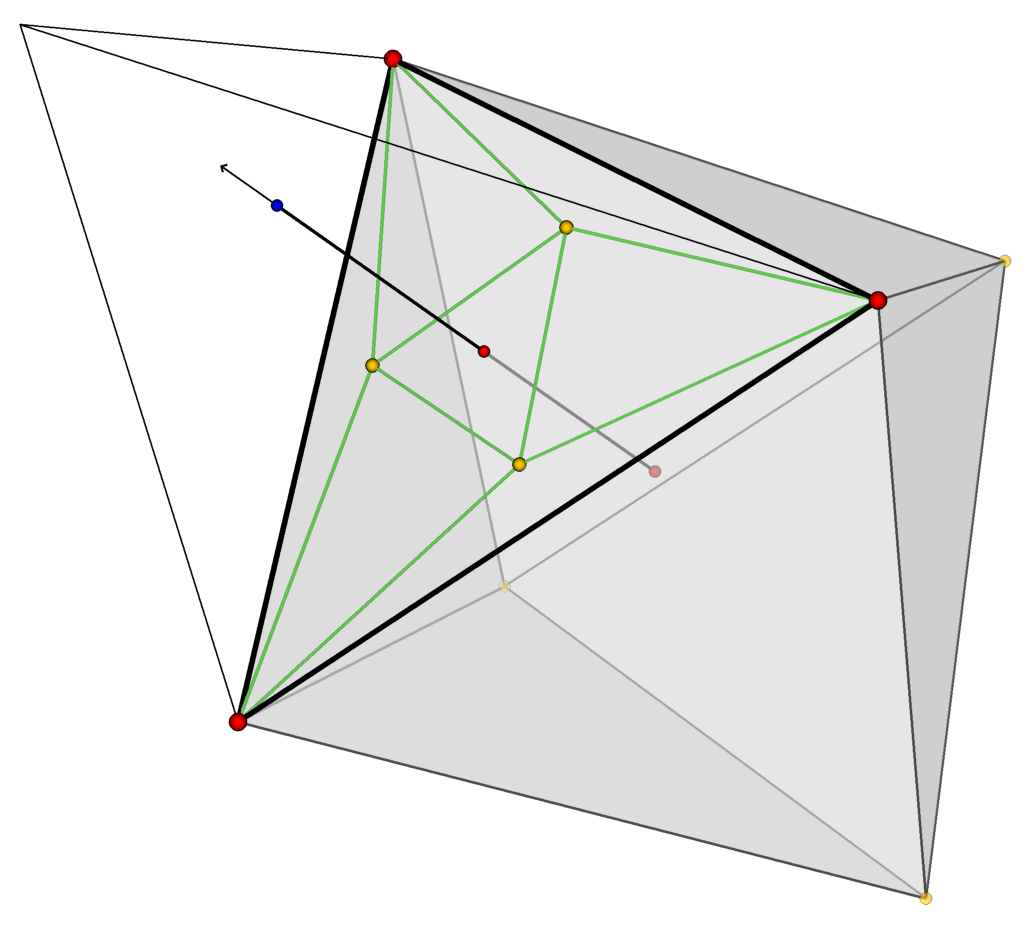} \put(28,70.5){$v$} \put(19.5,70){$r$}
\put(48,56.5){$w$}
  \end{overpic}
  \caption{Construction of a Schlegel diagram of the regular octahedron. The viewpoint $v$
lies in the region beyond the projection facet. It is described by its relative position
on the ray $w+\RR_{\ge 0}r$. }
  \label{fig:schlegel_viewpoint}
\end{figure}

It is more subtle to move the point beyond.  The simple reason for this is that it does
not appear in the $3$-dimensional Schlegel diagram nor its affine span.  Thus it must
be moved implicitly. In the beginning we choose a fixed point $w$ in the relative
interior of the projection facet $F$.  And we also choose a vector $r$ such that
$w+\RR_{\ge 0} r\cap\partial P=\{w\}$.  That is, $r$ points from $w$ towards the set of
points beyond~$F$. For $r$ we can take an outward pointing normal vector of~$F$; see
Figure~\ref{fig:schlegel_viewpoint}.

\begin{figure}[t]%[htbp]
  \centering
  \begin{overpic}[width=.48\linewidth]{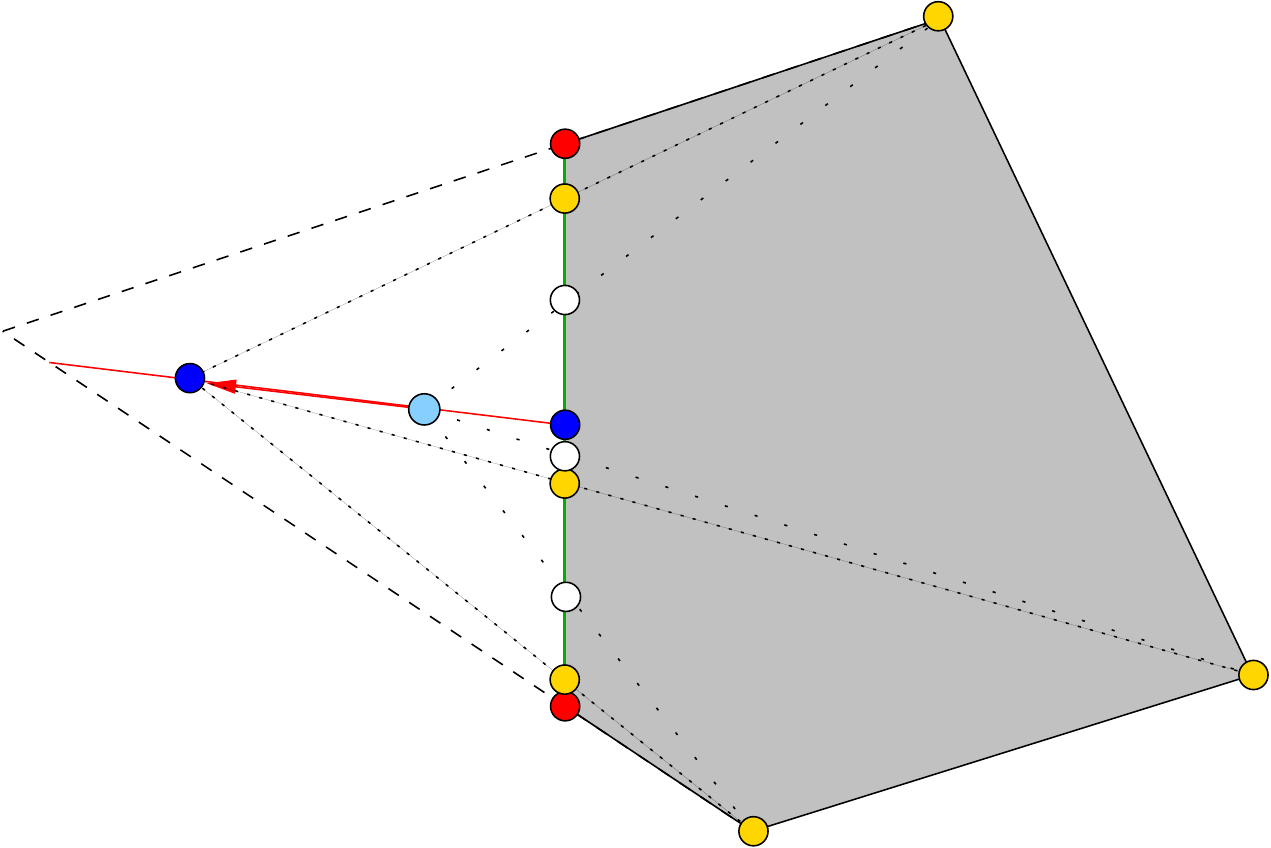} \put(7,39.5){$v'$} \put(27,37){$v$}
\put(47,33){$w$}
  \end{overpic}
  \caption{The construction of the Schlegel diagram of a pentagon (shaded) on one of its
facets (vertical edge).  The zoom parameter allows to move the old viewpoint $v$ to a new
viewpoint $v'$ on the ray connecting viewpoint and the point $w$ on the facet.}
  \label{fig:schlegel_zoom}
\end{figure}

Two cases are to be distinguished. First let us assume that the set of points beyond~$F$
is bounded, that is, there is a maximal $\lambda>0$ such that $w+\zeta\lambda r$ is
beyond~$F$ for all $\zeta\in(0,1)$.  We call $\zeta$ the \emph{zoom} value of the
viewpoint
$v=w+\zeta\lambda r$ with respect to $w$ and $r$.  The zoom value can be changed directly
in the graphical interface to obtain the Schlegel diagrams for viewpoints on the segment
$w+[0,\lambda]r$; see Figure~\ref{fig:schlegel_zoom}.  If, however, there is no such
maximal $\lambda$ the interval of~$\zeta$ is mapped from $(0,1)$ to $(0,+\infty)$ by
$\zeta \mapsto \tfrac{\zeta}{1-\zeta}$.

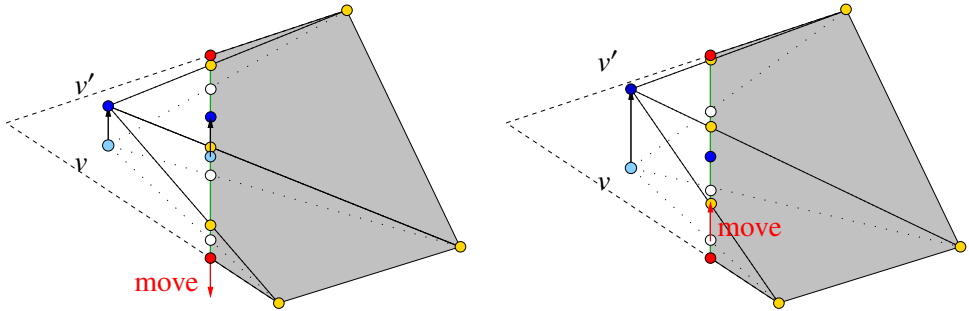
\begin{figure}[htbp] \centering
  \begin{minipage}{.48\linewidth}
    \begin{overpic}[width=\linewidth]{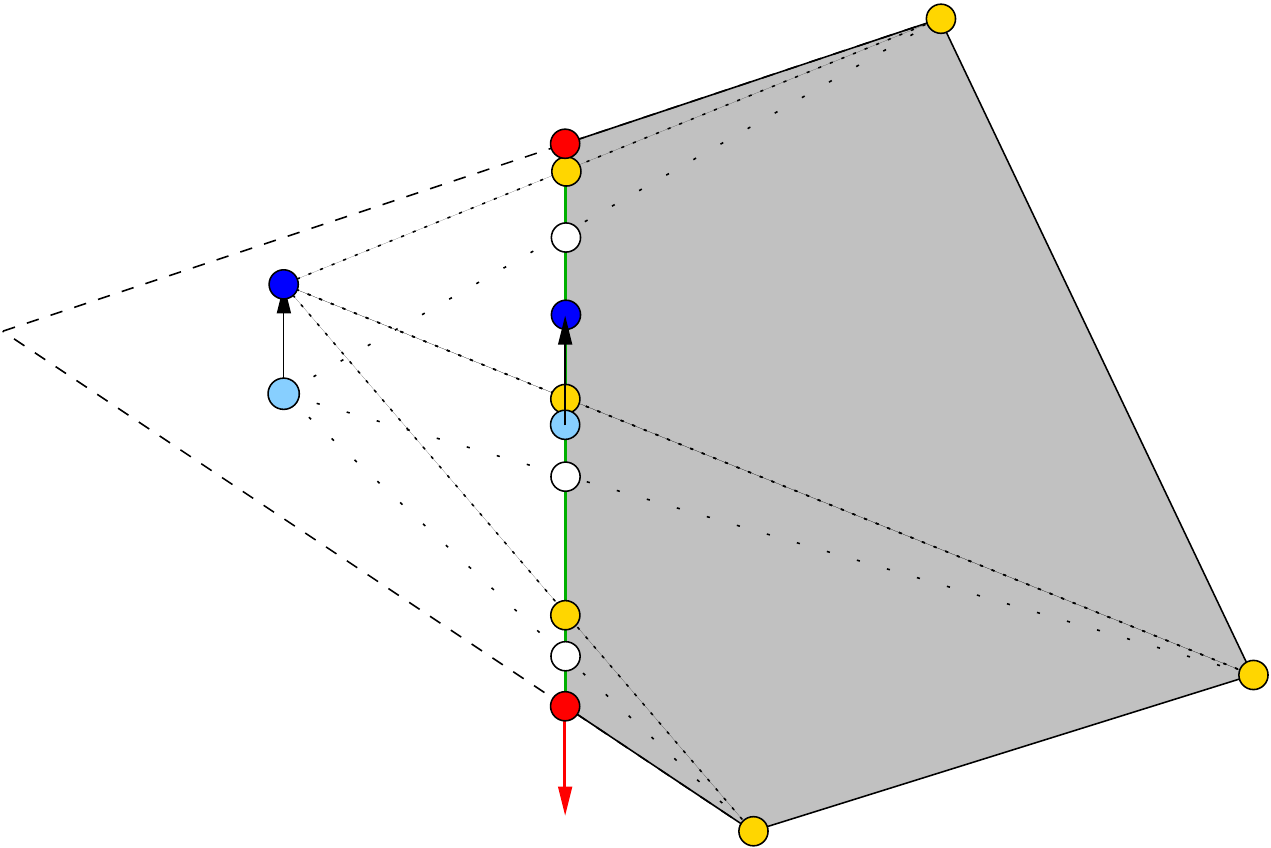} \put(15,30){$v$}
\put(15,47){$v'$} \put(28,4){\textcolor{red}{move}}
    \end{overpic}
  \end{minipage} \hfill
  \begin{minipage}{.48\linewidth}
    \begin{overpic}[width=\linewidth]{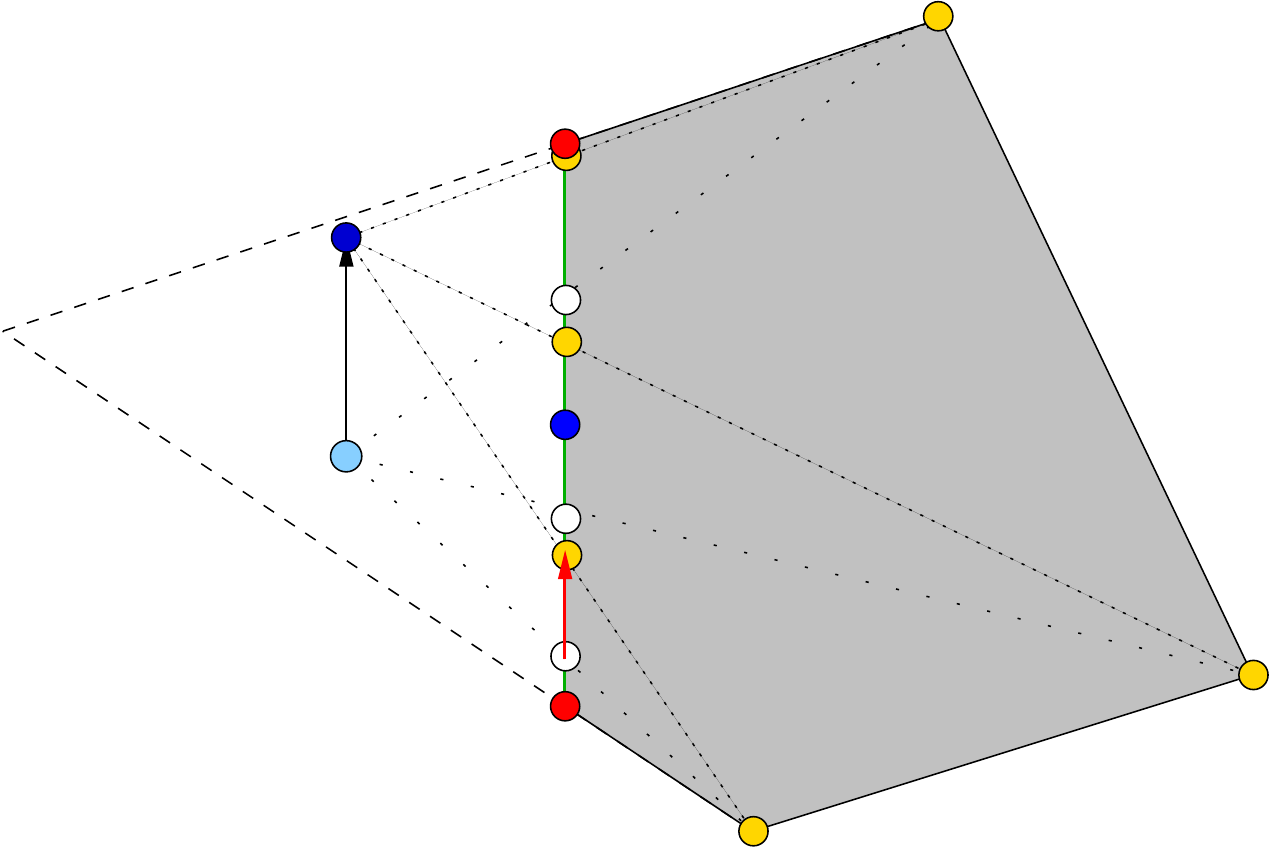}
      \put(20,26){$v$}
      \put(20,52){$v'$}
      \put(46,16){\textcolor{red}{move}}
    \end{overpic}
  \end{minipage}
  \caption{A pentagon with its Schlegel diagram on one of its facets.  The movement of the
vertices in the Schlegel diagram is translated into movement of the viewpoint, which does
not appear in the diagram. The movement depends on whether the point moved lies on the
projection facet (left) or not (right).}
  \label{fig:schlegel_move}
\end{figure}

Other Schlegel diagrams are obtained by dragging individual points in the projection.
There are surely different ways to move the viewpoint in the entire region beyond the
projection facet by interpreting the dragging of the vertices of the Schlegel diagram. Our
choice proved to be intuitive and very useful for our applications. It works as follows.
If a vertex of the projection facet $F$ is dragged, the current viewpoint~$v$ and the
point~$w$ are moved in opposite directions; see Figure~\ref{fig:schlegel_move}~(left).
If, however, a point $\pi(x)$ is moved, where $x$ is a vertex not belonging to~$F$, only
the viewpoint~$v$ is modified.  Then the dragged point becomes the projection of $x$ in
the modified Schlegel diagram; see Figure~\ref{fig:schlegel_move}~(right).  Since we
require that $v$ can only be moved parallel to $F$, this uniquely defines the new
viewpoint $v'$ as the intersection of the line through $x$ and $\pi(x)$ with the affine
hyperplane through~$v$ which is parallel to~$F$.

This approach allows to produce all possible Schlegel diagrams.  In the implementation it
is always verified that the new viewpoints lead to valid Schlegel diagrams, that is, the
new viewpoints remain to be points beyond $F$.

\section{Pseudo-Physical Models}\label{sec:pp_models}

\noindent So far we were concerned with the visualization of (low-dimensional) polytopes,
where we had natural ways of visualizing, either directly or via Schlegel diagrams.  In
higher dimensions or if only the face poset of a polytope is given (but no coordinates)
other techniques are required.  In the following we concentrate on visualizing the
vertex-edge graph of a polytope.  At least if the polytope is simple this can be expected
to be fruitful in view of Theorem~\ref{thm:BlindMani}.

A frequently followed approach is via models copied from physics. This way often nice
drawings can be obtained since the inherent symmetry properties of physical laws tend to
retain the abstract symmetry of a graph in its drawings.  Another great advantage of such
pseudo-physical models is that forces may be added to improve an existing model. On the
other hand, evolution in time of the pseudo-physical model has to be approximated and
convergence to a stable state is not guaranteed due to the discretization of time; for a
more thorough discussion see Fruchtermann and
Reingold~\cite{fruchtermann92:_graph_drawin_force_placem} or Tollis et
al.~\cite[Chapter~10]{TollisEtAl99}.

\subsection{Attracting and Repellent Forces}\label{subsec:attr_and_repell_forces}

\noindent We want to visualize a finite graph $G=(V,E)$.  The naive idea is to assign
random coordinates to each node $v\in V$ and then to let some ``forces'' act on the points
until an equilibrium is reached.  The \emph{neighbors} of a node~$v\in V$ form the set $
N(v) = \smallSetOf{w\in V}{\{v,w\}\in E} $, its \emph{closed neighborhood} is $ N[v] =
N(v)\cup\{v\}$.

At a given time each node of $G$ is represented by some point in $\RR^3$.  In the formulae
below we will identify each node with its coordinates in $3$-space.  Firstly we define a
repellent force
$\tfrac{-\delta_{\text{rep}}}{\|w-v\|^3}(w-v)$ pushing~$v$ away from any non-neighboring
node $w\in V(G)\setminus N[v]$.  One may think of this repellent force as resulting from a
(kind of) negative electronic charge carried by the vertices.  Then~$\delta_{\text{rep}}$
is the electrostatic constant in the repellent force.  On the other hand, the edge
$\{v,w\}$ for each adjacent node $w\in N(v)$ pulls~$v$ towards~$w$ like a stretched
spring.  Thus an attracting force of $\Big(\tfrac{1}{\ell_{\{v,w\}}} -
\tfrac{1}{\|w-v\|}\Big)(w-v)$ acts on~$v$, where $\ell_{\{v,w\}}$ is the \emph{desired
length} of the spring modeling the edge $\{v,w\}$. Depending on the situation this desired
length may be constant (for instance, if no coordinates are known) or it may reflect some
geometric properties (such as the Euclidean distance in some high-dimensional space).
Note that the exact formulation of the attracting and repelling forces described above
does not strictly reflect some physical model but also takes into account experimental
fine tuning. Summing up the attracting and repelling forces for~$v$ yields
\begin{equation}\label{eq:f_v} f_v \ = \ \sum_{w \in V\setminus N[v]}
\frac{-\delta_{\text{rep}}}{\|w-v\|^3}(w-v) \ + \ \sum_{w \in N(v)} \Big(\tfrac{1}{\ell_{\{v,w\}}} - \tfrac{1}{\|w-v\|}\Big) (w-v)\, .
\end{equation} This defines a discrete vector field on the finite point set representing
the node set~$V$.

In order to facilitate convergence we add inertia and viscosity to our dynamical system.
Adding inertia means to take the first derivative of the motion of each node into account.
Adding viscosity means to systematically decrease the total energy of the system.  This
will be achieved by scaling down the inertia by a constant~$\delta_{\text{visc}}$.  The
dynamics are modeled by a rather crude discretization of time.  This can result in
convergence problems, which may be avoided by using more advanced numerics. However, for
the scope of this paper a simple discretization will suffice.

Let~$v_i$ denote the coordinate vector of a node~$v$ at time~$i$. Then the new coordinate
vector~$v_{i+1}$ is given by
\[
v_{i+1} \ = \ v_i + f_{v_i} + \delta_{\text{visc}} (v_i - v_{i-1}) \, .
\]
As a start configuration choose any random distribution of the nodes on the unit sphere with inertia zero. The
simulation is run till the absolute fluctuation
\[
\max_{v \in V(G)} \|v_i - v_{i+1}\|^2
\]
drops below some fixed threshold.  So far the pseudo-physical model is quite standard, and it does not reflect any
special properties of polytopal graphs.

\begin{figure}[htbp]\centering
  \includegraphics[width=.29\textwidth]{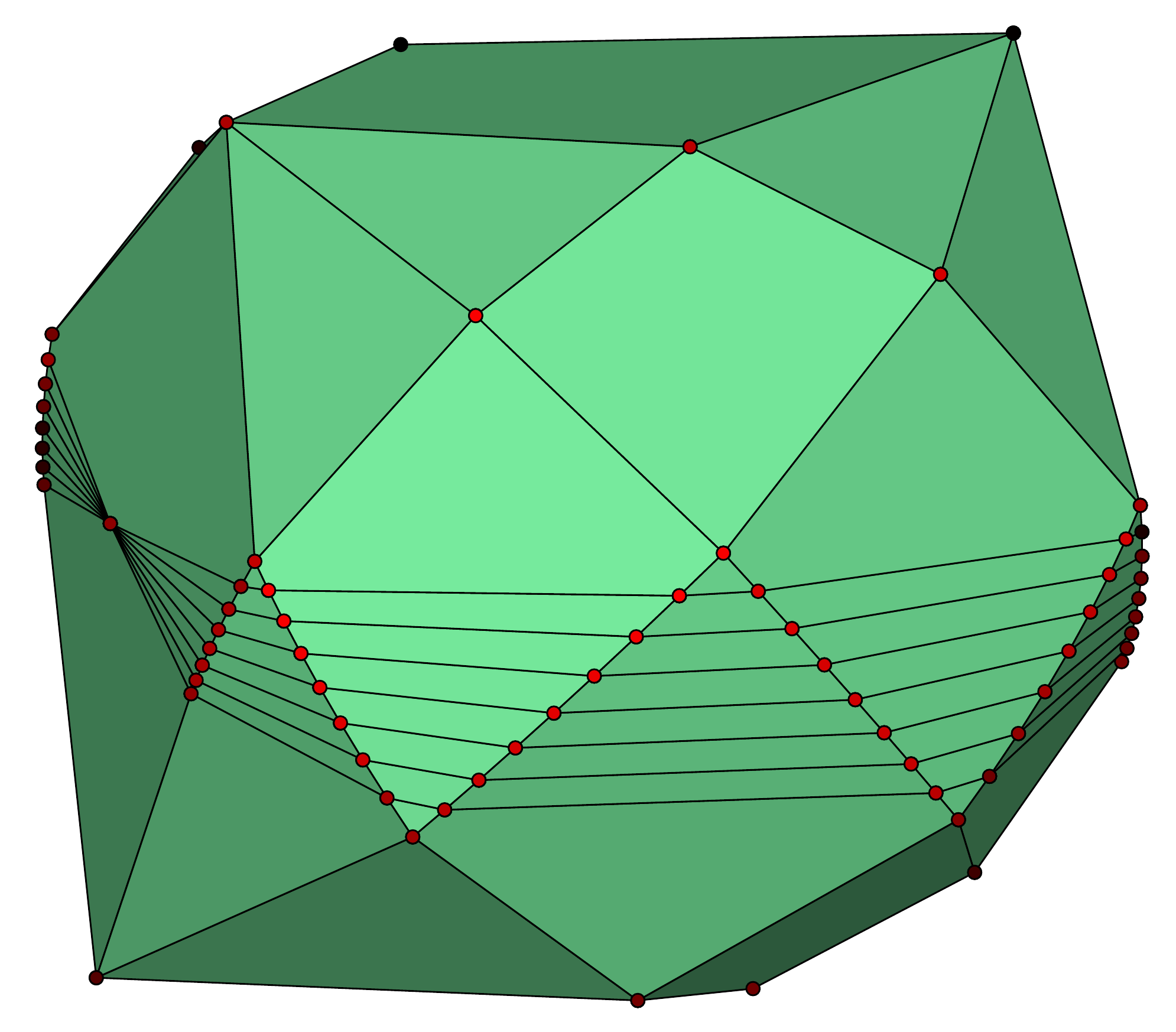}
  \hspace{.3cm}
  \raisebox{.15cm}{\includegraphics[width=.30\textwidth]{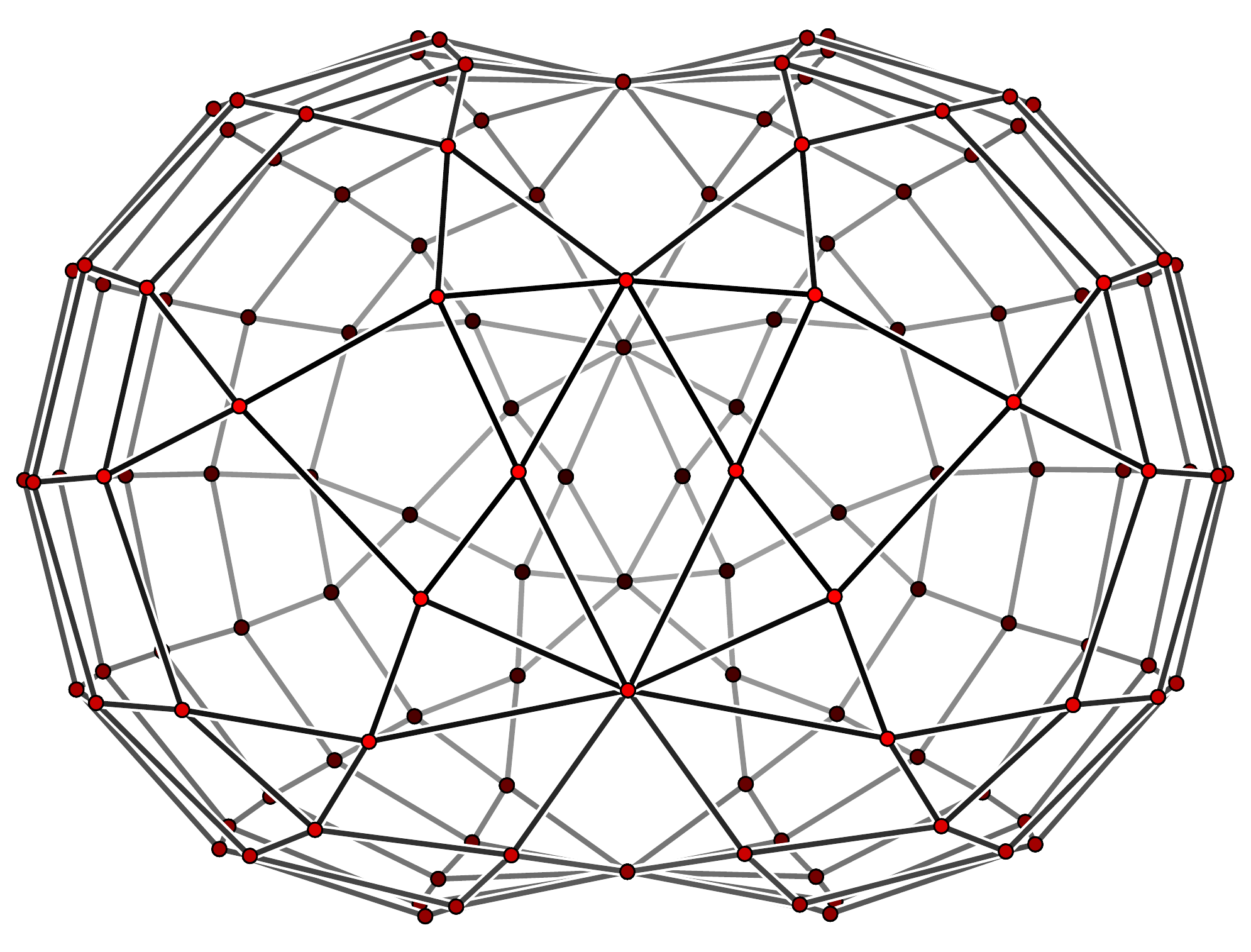}}
  \hspace{.3cm}
  \raisebox{.06cm}{\includegraphics[width=.33\textwidth]{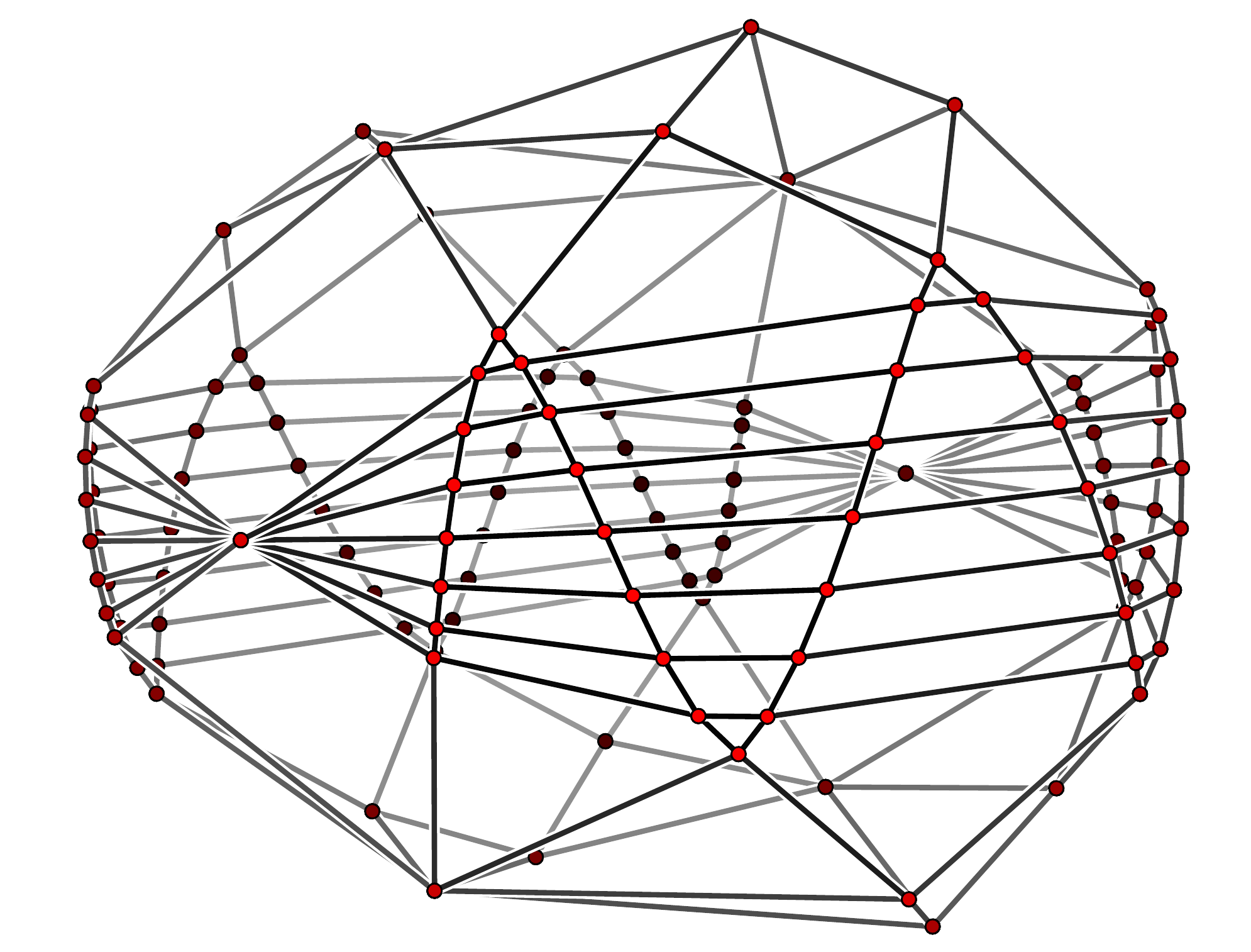}}
  \caption{A dual 3-zonotope embedded with its original coordinates (left) and two spring
    embeddings of its graph (center and right).  In the center drawing all edges are given
    the same desired edge lengths, on the right an attempt to embed each edge with its
    original length.}
  \label{fig:edge_length}
\end{figure}

\begin{ex}
  Bern, Eppstein et al.~\cite{BerEppGui-ESA-95,eppstein:_ukrain_easter_egg} construct
  3-dimensional zonotopes with central 2D-sections of quadratic size with respect to the
  number of zones (the data defining a zonotope). These $3$-zonotopes are remarkable since
  naively one might expect only a linear number of vertices in any section of a
  ``typical'' $3$-dimensional zonotope. Their duals are dual 3-zonotopes with~$n$ zones
  and a $2$-dimensional affine image (a \emph{$2$D-shadow}) with~$\Omega(n^2)$ vertices.
  Koltun~\cite[Problem~3]{OWReport} asked for a generalization, that is, dual
  $d$-zonotopes with~$n$ zones and 2D-shadows of size~$\Omega(n^{d-1})$ (for fixed~$d$).
  Note that there is a trivial upper bound of~$O(n^{d-1})$ for the size of a maximal
  2D-shadow since dual $d$-zonotopes correspond to affine $(d-1)$-dimensional hyperplane
  arrangements. Such dual $d$-zonotopes with 2D-shadow of size~$\Omega(n^{d-1})$ exist by
  \cite{roerig-witte-ziegler2007}.  A $3$-dimensional example is used in
  Figure~\ref{fig:edge_length} to illustrate the use of different desired edge lengths in
  the spring embedder.  This is a particularly challenging case for our pseudo-physical
  model because of the great length differences.
\end{ex}

However, one can modify such a pseudo-physical model by inventing further forces.  Here we
will pursue how to visualize the vertex-edge graph $G=\Gamma(P)$ of a polytope
$P\subset\RR^d$ if additionally a linear objective function $\lambda:\RR^d\to\RR$ is
given.  Without loss of generality we can assume that $\lambda$ projects a point
$x\in\RR^d$ to its last coordinate $x_d$.  The effect of the additional force $f_F$ may
be interpreted as a (vertical) linear field: Every vertex~$v\in G$ tries to adopt its
$x_3$-coordinate~$v_3$ (relative to the other vertices) according to its (relative) value
with respect to~$\lambda$.  To this end the center of gravity $\bar{v} =
\tfrac{1}{|V|}\sum_{v\in V}v$ of all vertices $V=V(G)$ and the average value
$\bar{\lambda}=\tfrac{1}{|V|}\sum_{v\in V}\lambda(v)$ is computed. Then the additional
vertical force ($e_3$ being the third unit vector)
\[
\big((\lambda(v)-\bar{\lambda})-(v-\bar{v})_3\big) e_3\, ,
\]
scaled by some constant $\delta_{\text{lin}}$, is added to the Equation~\eqref{eq:f_v}.

All the constants mentioned, such as $\delta_{\text{rep}}$, $\delta_{\text{visc}}$, and
$\delta_{\text{lin}}$, are non-negative.  They must often be chosen interactively in order
to balance the forces according to esthetic needs; this functionality is provided by
\polymake's interfaces to \JavaView and \JReality.

\begin{ex}
  The \emph{Klee-Minty cube}
  \cite{klee-minty1972,eg-models_goldfarb} is a $d$-dimensional polytope which is
  combinatorially isomorphic to the regular $d$-cube. It is defined as the set of
  admissible solutions of the linear program (maximizing~$x_d$) given by the~$2d$
  inequalities
  \begin{equation}\label{eq:goldfarb}
    \begin{aligned}
      0 &\leq \!\!\!&x_1 &\leq 1\\
      \tfrac{1}{3}x_i &\leq \!\!\!&x_{i+1} &\leq 1 -
      \tfrac{1}{3} x_i \quad \text{for $1\leq i<d$.}
    \end{aligned}
  \end{equation}
  The Klee-Minty cube has an \emph{ascending path} of length~$2^d$, that is, there exists
  a directed path of length~$2^d$ in its graph such that any vertex of this path has
  greater $x_d$-value than its predecessor.  This provides an example of a polytope with
  an exponentially long (with respect to the number of defining halfspaces) ascending path
  and thus a ``bad case'' for the simplex algorithm; see Figure~\ref{fig:gf3} (left).

  \begin{figure}[htbp]\centering
    \includegraphics[width=.32\textwidth]{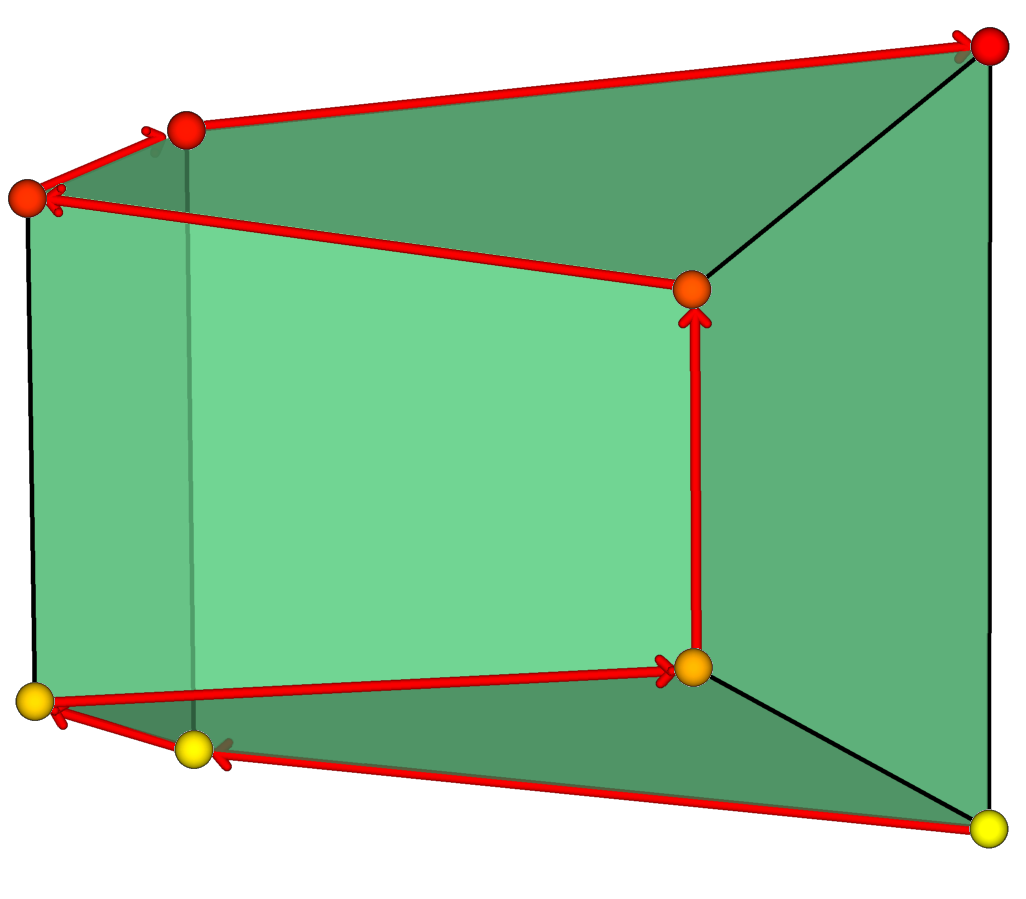}
    \includegraphics[width=.32\textwidth]{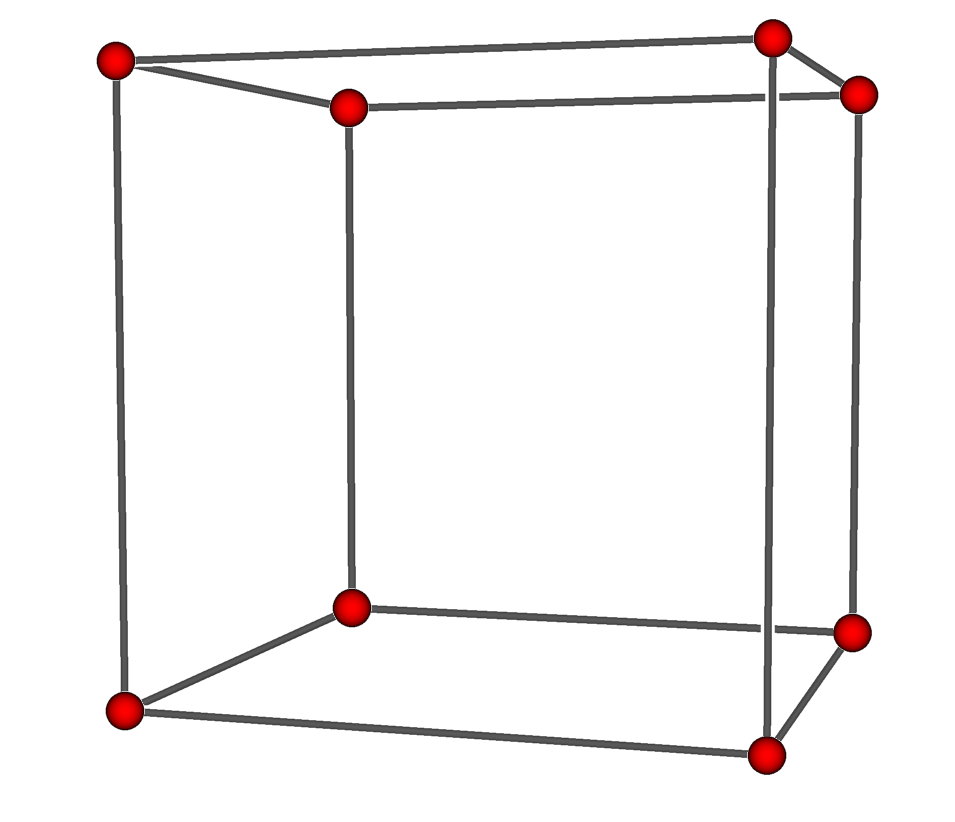}
    \includegraphics[width=.32\textwidth]{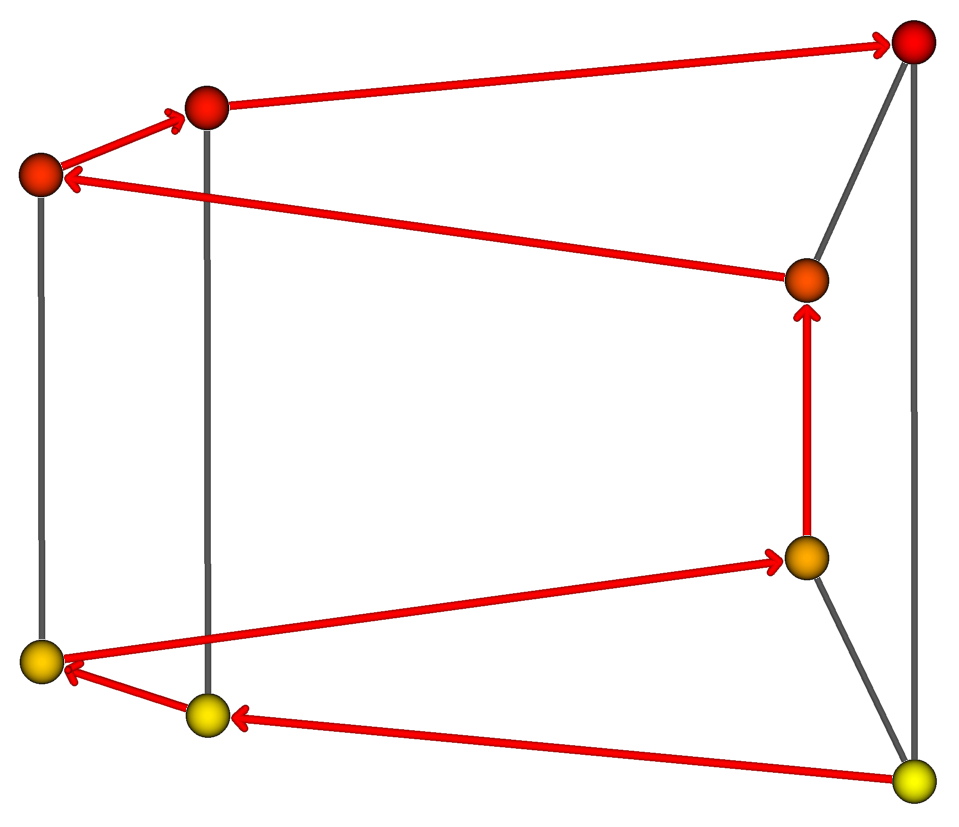}
    \caption{The Klee-Minty 3-cube (with its realization given by
      Equation~\eqref{eq:goldfarb}) and an ascending Hamiltonian path on the left. On the
      right two embedding of its graph, first without additional forces, then with an
      additional vertical force determined by the linear objective function
      $\RR^3\to\RR:x\mapsto x_3$.
      \label{fig:gf3}}
  \end{figure}
  
  We choose the 3-dimensional Klee-Minty cube to illustrate the effect of the vertical
  force governed by the linear objective function $\RR^3\to\RR:x\mapsto x_3$. Of course, a
  3-dimensional polytope may be visualized directly; see Figure~\ref{fig:gf3} (left).
  However, Figure~\ref{fig:gf3} (center and right) exhibits the effect of the additional
  vertical force: The drawing in the center does not reflect the particular realization
  of the Klee-Minty cube, the embedding on the right on the other hand clearly shows the
  ascending Hamiltonian path.
\end{ex}

\begin{figure}[htbp]\centering
  \includegraphics[width=.42\textwidth]{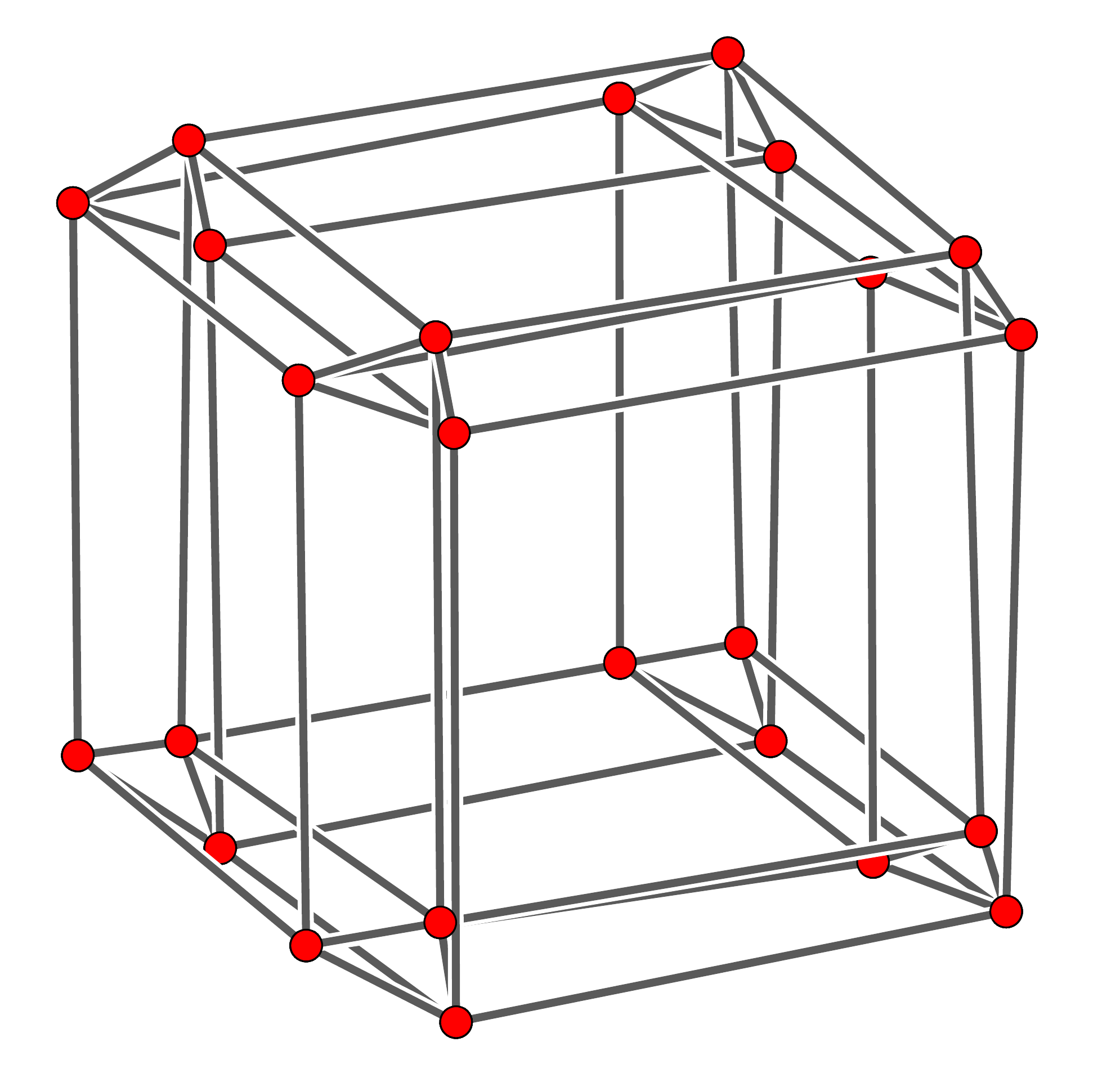}\qquad
  \raisebox{.5cm}{\includegraphics[width=.48\textwidth]{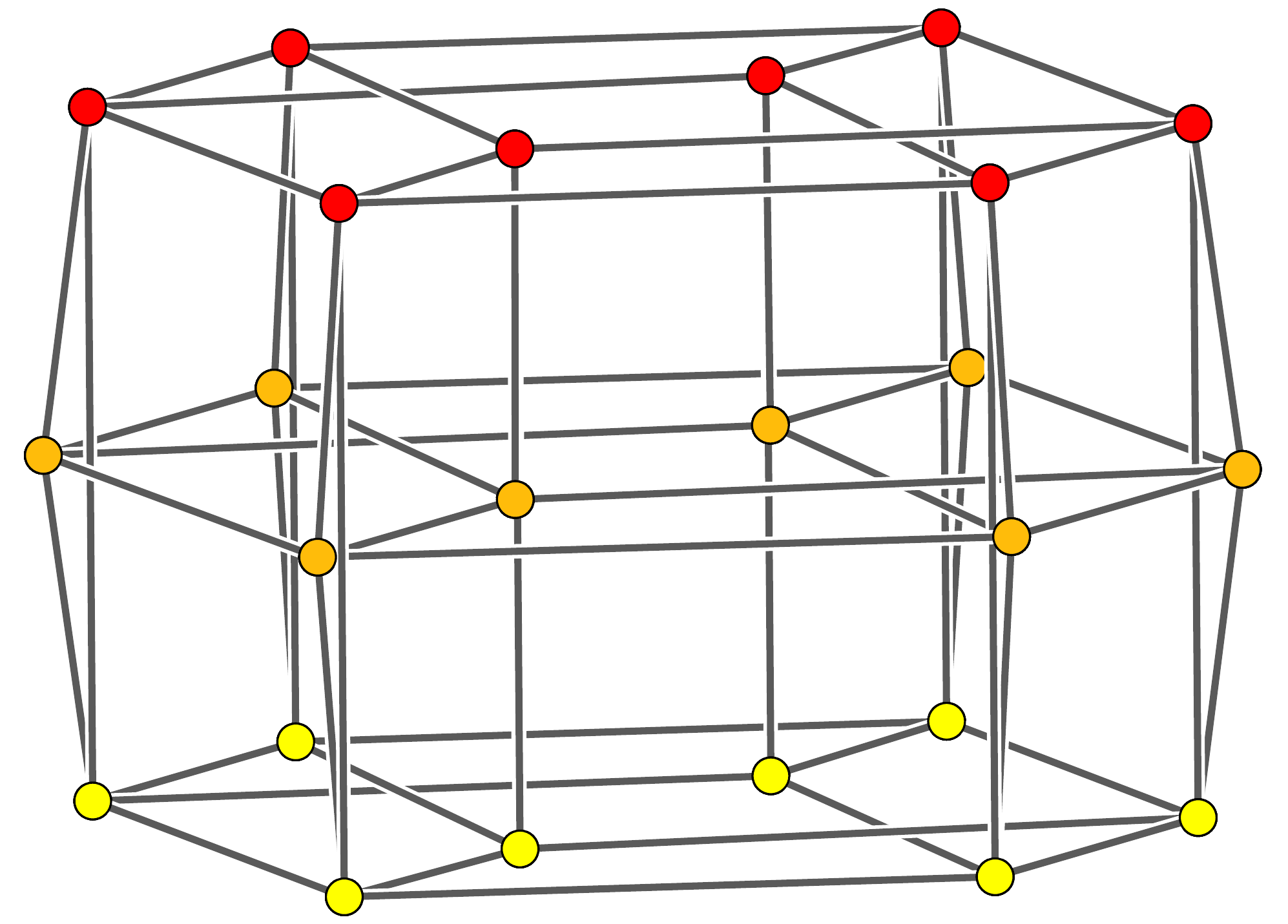}}
  \caption{Two embeddings of the graph of the product of a triangle and the 3-cube. The
    product $\Delta_2\times C_3$ has six ``triangle~$\times$ square'' facets (emphasized on
    the left) and three ``edge~$\times$ cube'' facets (emphasized on the right).}
  \label{fig:s2xc3}
\end{figure}

\begin{ex} We want to visualize the product $\Delta_2\times C_3$ of a triangle and the
  $3$-dimensional unit-cube.  This is a $5$-dimensional simple polytope.  We choose a
  linear objective function~$\lambda:\RR^5\to\RR$ such that~$\lambda$ evaluates to
  distinct values on the vertices of $\Delta_2$, but does not distinguish in between the
  vertices of~$C_3$.  For example let $\lambda(x)=l_1x_1+l_2x_2$ for random
  values~$l_1,l_2\in\RR$. Hence for a vertex~$v$ of~$\Delta_2$ all vertices of the type
  ``$v$ times any vertex of~$C_3$'' are embedded with the same $x_3$-coordinate by the
  vertical force.  Figure~\ref{fig:s2xc3} depicts a spring embedding of the graph of
  $\Delta_2\times C_3$ without additional forces on the left, and on the right
  with~$\lambda$ taken into account. The image on the right clearly shows the three (flat
  horizontal) copies of~$C_3$ corresponding to the three vertices of~$\Delta_2$, which in
  turn provide the second factor of the three ``edge times cube'' facets
  of~\mbox{$\Delta_2\times C_3$}.
\end{ex}

\subsection{Rubber Bands}\label{sec:rubber_bands}

\noindent In the following we take a different approach to embedding of graphs
in~$\RR^3$. Rather than working with a dynamic model as in
Section~\ref{subsec:attr_and_repell_forces}, we present a classical method due to
Maxwell~\cite{maxwell1864:on_reciprocal_figures} and Cremona~\cite{Cremona1890}.  Here a
graph~$G$ is embedded into~$\RR^3$ by solving a system of linear equations. This method
requires some nodes~$\Phi$ already embedded in~$\RR^3$.  The graph~$G$ will be embedded in
the affine subspace spanned by the nodes in~$\Phi$. Hence it is often useful to require
that~$\Phi$ contains at least four nodes which span the whole space.  Further we
assume~$G$ to be connected, since we may embed different connected components
individually.

We picture an edge~$e\in E(G)$ as a spring, or a rubber band, of length zero (if not
stretched); the edge $e$ has an individual spring constant~$\delta_e$. After fixing
coordinates for the nodes in~$\Phi$ we let the rubber bands pull the remaining nodes to an
equilibrium. This equilibrium is attained by minimizing the total energy~$E$ of the system
of rubber bands.  To this end let~$v_1$,~$v_2$, and~$v_3$ be the coordinates in~$\RR^3$ of
a node~$v\in V(G)$ and the energy of a rubber band representing an embedded edge~$e$
is~$\delta_e\|e\|^2$. Hence we have
\[
E \ = \ \sum_{e \in E(G)} \delta_e \|e\|^2 \ = \ \frac{1}{2} \sum_{v\in V(E)} \sum_{w\in N(v)}
\delta_{\{v,w\}} \left\|
\begin{pmatrix}v_1-w_1\\v_2-w_2\\v_3-w_3\end{pmatrix}
\right\|^2 \, .
\]
The total energy~$E$ is a quadratic function in $3 |V(G) \setminus \Phi|$ variables
$v_1$, $v_2$, and $v_3$ for all $v \in V(G) \setminus \Phi$ and thus~$E$ has a unique
minimum. Partial differentiation with respect to $v_1$ yields
\[ \frac{\partial E}{\partial v_1} \ = \ 2 \sum_{w \in N(v)} \delta_{\{v,w\}} (v_1 - w_1)
\, ,
\] and likewise for $v_2$ and $v_3$. Requiring equilibrium, that is, $\frac{\partial
E}{\partial v_1} = \frac{\partial E}{\partial v_2} = \frac{\partial E}{\partial v_3} = 0$
for all $v \in V(G) \setminus \Phi$, amounts to solving a system of $3 |V(G) \setminus
\Phi|$ linear equations to determine the values of $v_1$,~$v_2$, and~$v_3$ for all $v\in
V(G)\setminus \Phi$.

This technique was applied by Tutte~\cite{tutte63:how_to_draw_a_graph} to construct
crossing free embeddings of connected planar graphs into~$\RR^2$ with straight edges.
Here~$v_3$ is set to zero for all nodes~$v$ to obtain an embedding in~$\RR^2$.

\begin{figure}[t]\centering
  \includegraphics[width=.7\textwidth]{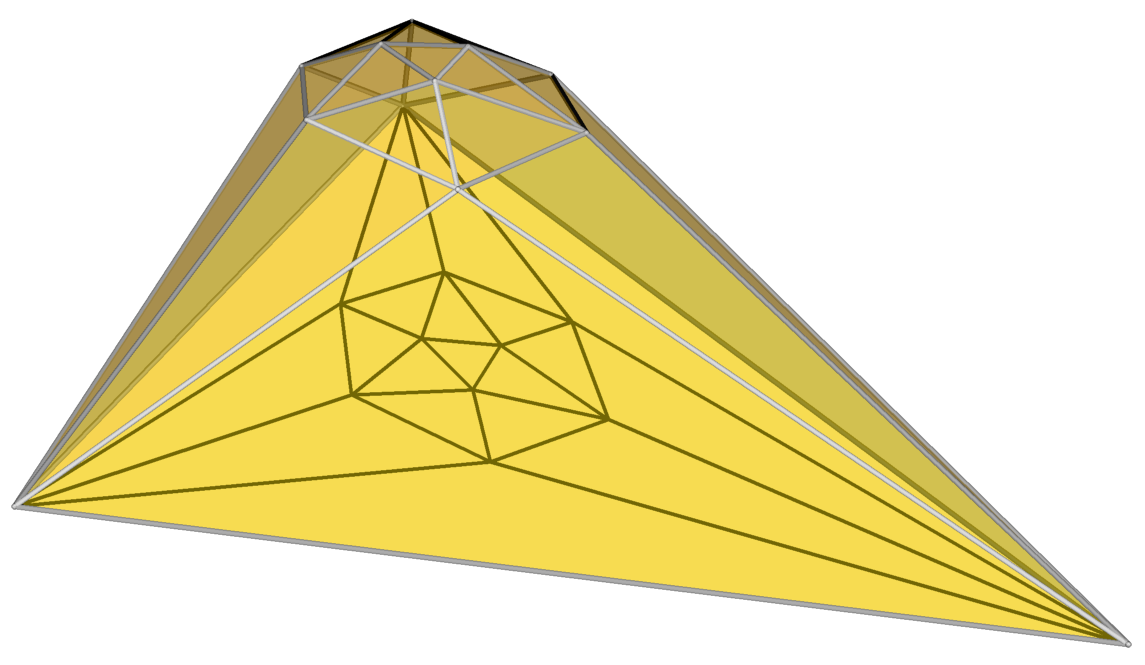}
  \caption{Non-regular realization of the icosahedron obtained from its graph via the
    rubber band method. The rubber band embedding of the graph is visualized as the shadow
    of the edges of the icosahedron.}
  \label{fig:tutte}
\end{figure}

The same method can also be used to prove the difficult direction of Steinitz' Theorem
(Theorem~\ref{thm:Steinitz}): The construction of a 3-polytope from a given planar,
3-connected graph. In the following we sketch the idea. Let $G$ be a $3$-connected planar
graph.  Then $G$ or its dual graph $G^*$ possesses a triangular face.  This can be derived
from Euler's theorem and double counting.  It suffices to prove that $G$ or $G^*$ is the
vertex-edge graph of some $3$-polytope $P$ because the polar dual of $P$ has the dual graph
as its vertex-edge graph.  So we can safely assume that $G$ has a triangular face.  Fix
its nodes in general position in~$\RR^2$. Now we embed~$G$ in~$\RR^2$ using Tutte's rubber
band method. According to Maxwell~\cite{maxwell1864:on_reciprocal_figures} such an
embedding may be lifted into~$\RR^3$, that is, there exists a convex function~$\RR^2
\rightarrow \RR$ which is linear on the faces of the planar embedding of~$G$. The
polytope~$P$ with~$G$ at its graph is the convex hull of the lifted nodes of~$G$. Since we
fixed the coordinates of a triangular face of~$G$ to begin with our rubber band method, no
new edges arise in the convex hull~$P$ and $G=\Gamma(P)$ holds. See Figure~\ref{fig:tutte}
for a realization of the icosahedron obtained from its graph via the rubber band
method. For a detailed proof see Richter-Gebert~\cite[Sect.\
13.1]{richter-gebert96:realization_spaces_of_polytopes}.

Notice that Tutte's rubber band embedding of $G$ is a Schlegel diagram of the lifted
polytope $P$ (with a viewpoint at infinity).

\section{Applications}

\noindent In contrast to what we discussed so far we will now study the visualization of
geometric objects which are more loosely connected to polytopes.  In most cases we will
deal with visualizing graphs by the pseudo-physical approach from
Section~\ref{subsec:attr_and_repell_forces} with additional forces which are specific to
the
application intended.

\subsection{Tight Spans of Finite Metric Spaces}
\label{sec:tight_spans}

\noindent Let $\delta:T\times T\to\RR_{\ge0}$ be a metric on a finite
set~$T=\{t_1,\dots,t_n\}$ of \emph{taxa}.  Then one can associate to it the convex
polyhedron
\[ P_\delta \ = \ \SetOf{x\in\RR^n}{x_i+x_j\ge \delta(t_i,t_j) \text{ for all $i,j$}} \, .
\] Because of $\delta(t_i,t_i)=0$ the polyhedron $P_\delta$ is contained in the positive
orthant $\RR_{\ge0}^n$, and hence it is \emph{pointed}, that is, it does not contain any
affine line.  Moreover, $P_\delta$ is always non-empty and unbounded, since the ray
$\RR_{\ge0}(M,\dots,M)$ is contained in~$P_\delta$, where~$M$ is the maximal value
attained by the metric~$\delta$.  The polytopal subcomplex of all those faces which are
bounded, the \emph{bounded subcomplex} of~$P_\delta$, is denoted by~$\cT(P_\delta)$.
Bandelt and Dress~\cite{MR858908} introduced the name \emph{tight span} for these objects,
but they already showed up earlier as the \emph{injective envelope} of a metric space in
the work of Isbell~\cite{MR0182949}.  Note that the bounded subcomplex of an unbounded
polyhedron is always contractible.  Tight spans of finite metric spaces are dual to
regular subdivisions of second hypersimplices.

The interest in this construction comes from the following simple observation.

\begin{prop}\label{prop:tree-like} The polytopal complex $\cT(P_\delta)$ is
$1$-dimensional, that is, it consists of edges only, if and only if the metric~$d$ is
tree-like.
\end{prop}

Here a metric is called \emph{tree-like} if it arises from a finite tree, that is, a
connected graph without cycles, with non-negative weights associated to the edges.
Between any two nodes in a tree there is a unique shortest path, and hence, by adding up
all the weights on a shortest path, this gives a distance function on the set of nodes.
As there are no cycles, there are no proper triangles, and the triangle inequality is
trivially satisfied.  The tight span of such a tree-like metric itself (essentially) is
the tree.  Moreover, the whole construction of the tight span depends on the metric in a
continuous way.  Hence, a tight span is a geometric object which can be used to measure in
how far a given metric is tree-like.  This very property can be exploited for an
interesting application.

The goal of phylogenetics, as a subject in biology, is to determine evolutionary kinship
among species or individual organisms.  Methods include the inspection of fossil samples,
the morphological analysis of extinct and existing species, and the study of ontogenetic
development on organisms.  One approach, which became feasible with the advent of modern
sequencing techniques, is to extract genetic distance information from alignments of DNA
or amino acid sequences.  Again there is a choice of mathematical models which try to
associate a tree with the given sequences.  But the beauty of the tight span approach lies
in the possibility to detect whether it makes sense to associate any tree with the given
metric to begin with.  In this sense tight span based techniques are less biased than
several other methods.

As far as the visualization is concerned several choices can be made.  For instance, one
can use the plain pseudo-physical model from Section~\ref{subsec:attr_and_repell_forces}.
However, this \emph{combinatorial visualization} does not provide us with images which are
meaningful in the context of phylogenetics.  For this it is better to set the desired edge
lengths in \eqref{eq:f_v} to edge lengths related to the distance function.  It turns out
that the taxa arise as specific vertices of the polyhedron $P_\delta$, and the proper
desired edge length is the distance with respect to maximum norm between two vertices of
$P_\delta$.  We call this the \emph{approximate metric visualization} of a tight span.

\begin{figure}[bth]\centering
  \raisebox{.3cm}{\includegraphics[width=.53\textwidth]{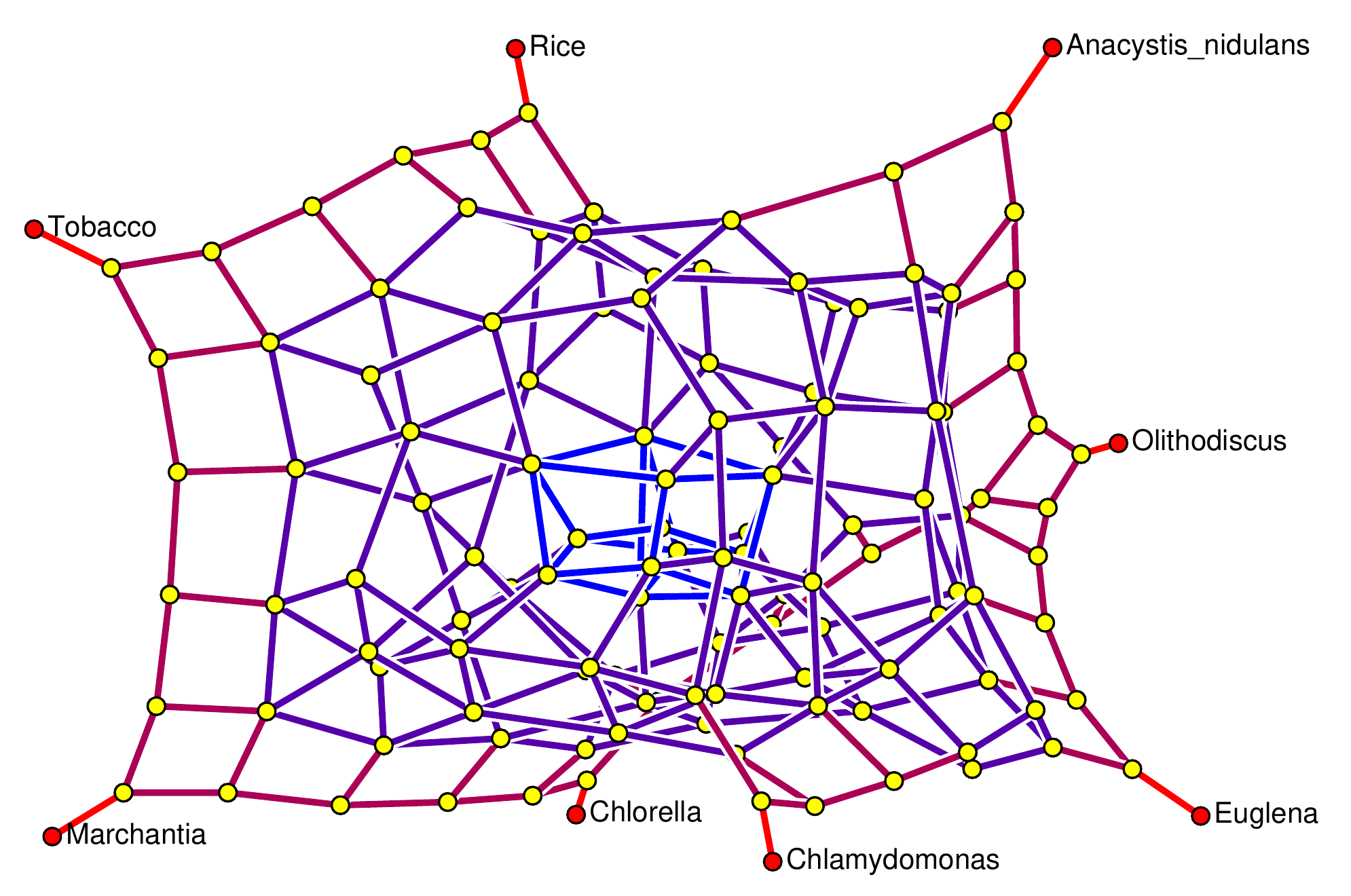}}\hfill
  \includegraphics[width=.47\textwidth]{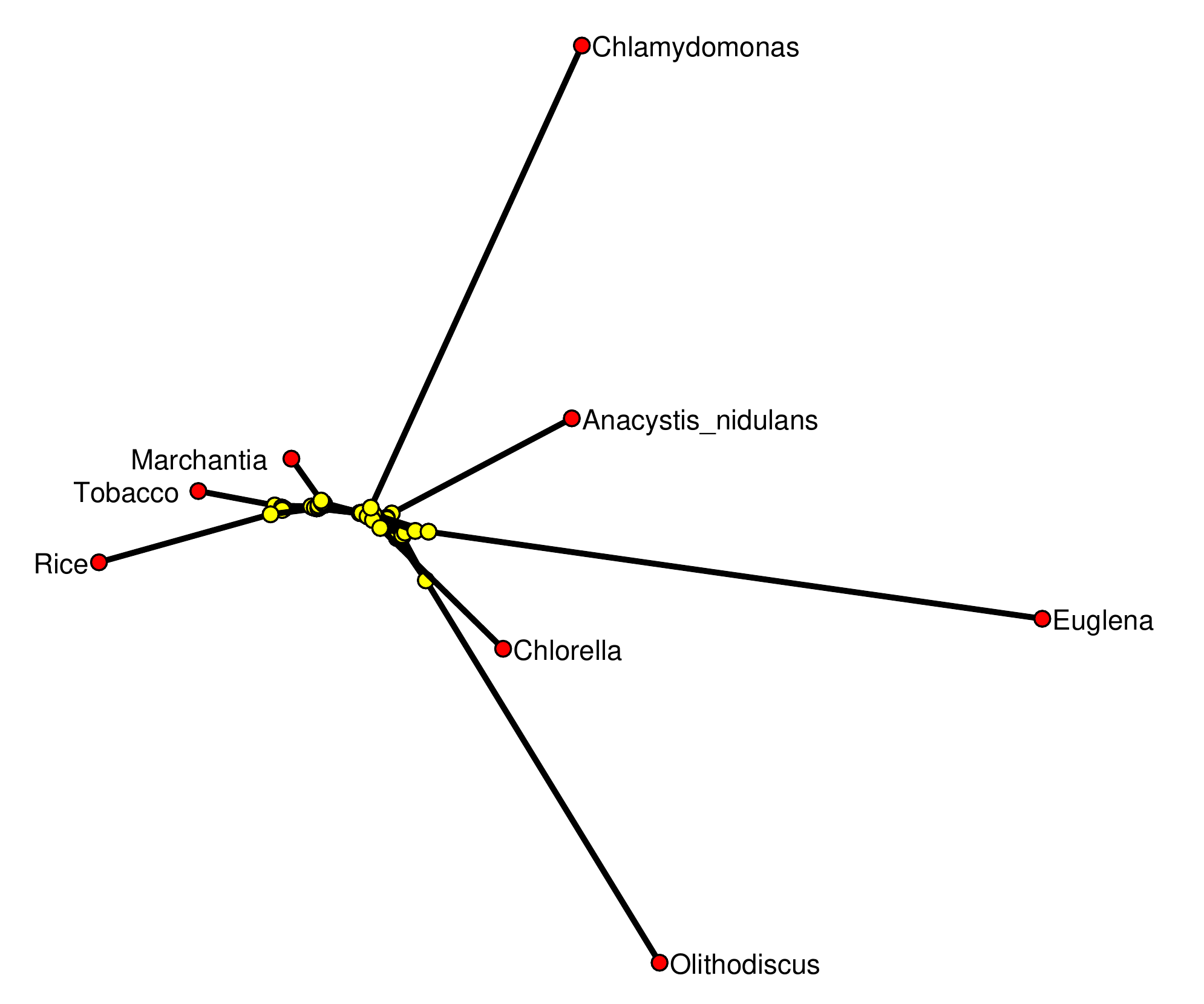}
  \caption{Left: Combinatorial visualizing of the tight span of five species of algae
    (chlamydomonas, chlo\-rel\-la, euglena, anacystis nidulans, olitho\-dis\-cus) and
    three non-algae (tobacco, rice, marchantia).  Right: Approximate metric visualization
    of the same tight span.}
  \label{fig:tightspan}
\end{figure}

\begin{ex} Consider a metric~$\delta$ on eight taxa which happens to be induced by aligned
RNA-samples of eight different species, five of which being algae.  The complete data set
is taken from the example file \texttt{algae.nex} which comes with
\texttt{SplitsTree}~\cite{splitstree}.  Note that the RNA samples used are very short:
They come from 920 bases each.  There is more than one way to compute a distance function
from these samples; here we used \texttt{SplitsTree}'s method \texttt{UncorrectedP} to
arrive at a metric given by the upper triangular matrix
\[
\small \left(
  \begin{array}{llllllll}
    0.0 & 0.026 & 0.029 & 0.112 & 0.078 & 0.136 & 0.123 & 0.141 \\
    & 0.0 & 0.041 & 0.121 & 0.088 & 0.144 & 0.132 & 0.145 \\
    & & 0.0 & 0.099 & 0.064 & 0.123 & 0.121 & 0.133 \\
    & & & 0.0 & 0.1 & 0.142 & 0.143 & 0.156 \\
    & & & & 0.0 & 0.116 & 0.118 & 0.116 \\
    & & & & & 0.0 & 0.159 & 0.135 \\
    & & & & & & 0.0 & 0.136 \\
    & & & & & & & 0.0 \\
  \end{array} \right) \, ,
\]
where the rows and columns are labeled with the following ordered list of taxa: tobacco,
rice, marchantia, chlamydomonas, chlo\-rel\-la, euglena, anacystis nidulans,
olitho\-dis\-cus.

  The left picture in Figure \ref{fig:tightspan} shows a combinatorial visualization of
the graph of $\cT(P_\delta)$.  The eight taxa can be associated with vertices of
$P_\delta$ and hence occur as nodes of the graph.  The colors of the edges visualize the
highest dimension of a bounded face containing that edge.  By default ``red'' refers to
dimension $1$ and ``blue'' to the maximum dimension $\dim\cT(P_\delta)$, which equals $4$
in this case.  The two shades of purple refer to dimensions $2$ and $3$, respectively.

  The picture to the right is the approximate metric visualization of $\cT(P_\delta)$.  No
color coding for the edges in this case.  The sample is too small to deduce much about the
evolutionary relationship between the five kinds of algae, but it already suffices to tell
the algae apart from the non-algae species.
\end{ex}

The visualization of tight spans used by Sturmfels and Yu is the combinatorial one
\cite{MR2097310}.

\subsection{Tropical Polytopes}

\noindent Consider a matrix $C=(c_{ij})\in\RR^{m\times n}$, and let
$W=\RR^{m+n}/\RR(1,1,\dots,1,-1,-1,\dots,-1)$ be a (quotient) vector space of dimension
$m+n-1$.  Each point in $W$ is written as a pair $(y,z)$ where $y$ takes the first $m$
coordinates, and $z$ the remaining $n$. Note that, in the quotient $W$, the equation
$(y+(1,1,\dots,1),z-(1,1,\dots,1))=(y,z)$ holds. Now
\begin{equation}\label{eq:tropical} T_C \ = \ \SetOf{(y,z)\in W}{y_i+z_j\le c_{ij}}
\end{equation} is a pointed unbounded convex polyhedron.  The bounded subcomplex 
$\cT(T_C)\subset W$ is the \emph{tropical polytope} generated by (the rows of) $C$.  The
vertices of $T_C$ are called the \emph{tropical pseudo-vertices} of $\cT(T_C)$ with
respect to the rows of $C$.  Among these there is a
unique inclusion-minimal subset which generates the tropical polytope: the
\emph{tropical vertices} of $\cT(T_C)$.  Tropical polytopes are dual to regular
subdivisions of products of simplices.  For details on the subject see
\cite{DevelinSturmfels04,Joswig05,BlockYu06}.

As it turns out, projecting $T_C$ onto the first $m$ coordinates (or onto the last $n$
coordinates) and clearing the first remaining coordinate (by adding a suitable multiple of
$(1,1,\dots,1)$) yields an affine isomorphism.  This way, one has a direct visualization
in $\RR^3$ if $\min(m,n)\le 4$.

\begin{figure}
  \begin{overpic}[width=.45\textwidth]{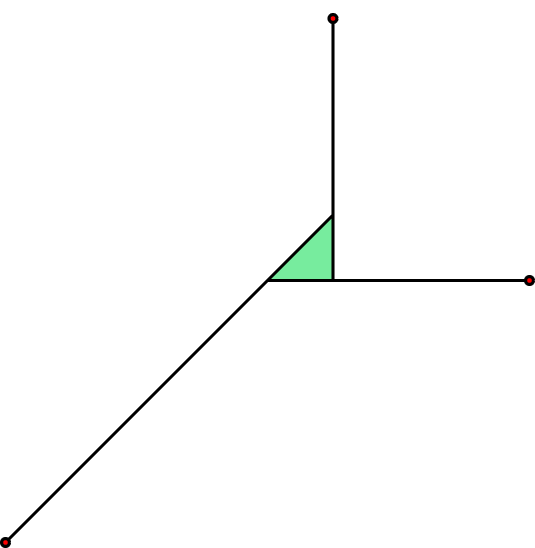} \put(4,-1){\small$(1,0,0)=(0,-1,-1)$}
\put(87.7,42.1){\small$(0,1,0)$} \put(63,95){\small$(0,\tfrac{1}{4},1)$}
  \end{overpic} \hfill
  \includegraphics[width=.45\textwidth,clip=true]{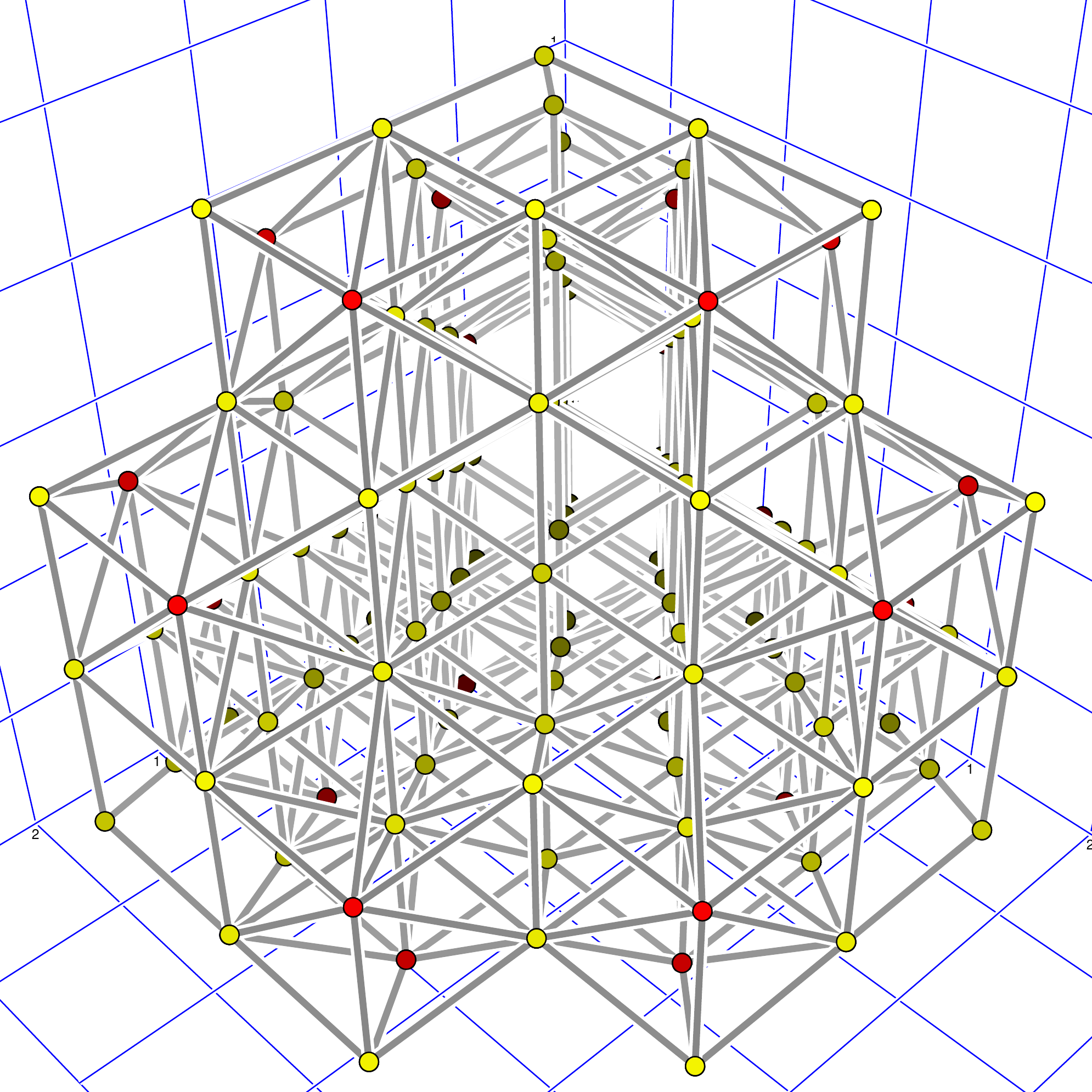}
  \caption{Left: Tropical triangle $\cT(T_C)$, where $C$ is defined in \eqref{eq:C}.
    Right: Tropical $3$-permutohedron.}
  \label{fig:trop-perm}
\end{figure}

\begin{ex} Letting
  \begin{equation}\label{eq:C} C \ = \
    \begin{pmatrix} 1 & 0 & 0 \\ 0 & 1 & 0 \\ 0 & \tfrac{1}{4} & 1
    \end{pmatrix}
  \end{equation} gives the tropical triangle $\cT(T_C)$, where $m=n=3$, in Figure
  \ref{fig:trop-perm} (left).  Its vertices correspond to the rows of the matrix $C$.
\end{ex}

\begin{ex} If we form the $(n!\times n)$-matrix whose rows are the permutation vectors
  from Example~\ref{ex:permutohedron} the construction \eqref{eq:tropical} gives the
  \emph{tropical $(n-1)$-permutohedron}. Figure \ref{fig:trop-perm} (right) shows the
  tropical $3$-permutohedron; the unit grid in the background is meant to provide a better
  idea about the spatial proportions.
\end{ex}

Since tropical polytopes are so similar to tight spans of finite metric spaces the same
visualization techniques can be applied in higher dimensions.  However, for the tropical
polytopes the combinatorial information is usually of interest.  This is why the
combinatorial visualization (with constant desired edge length) is preferred.

\begin{figure}
  \includegraphics[width=.8\textwidth]{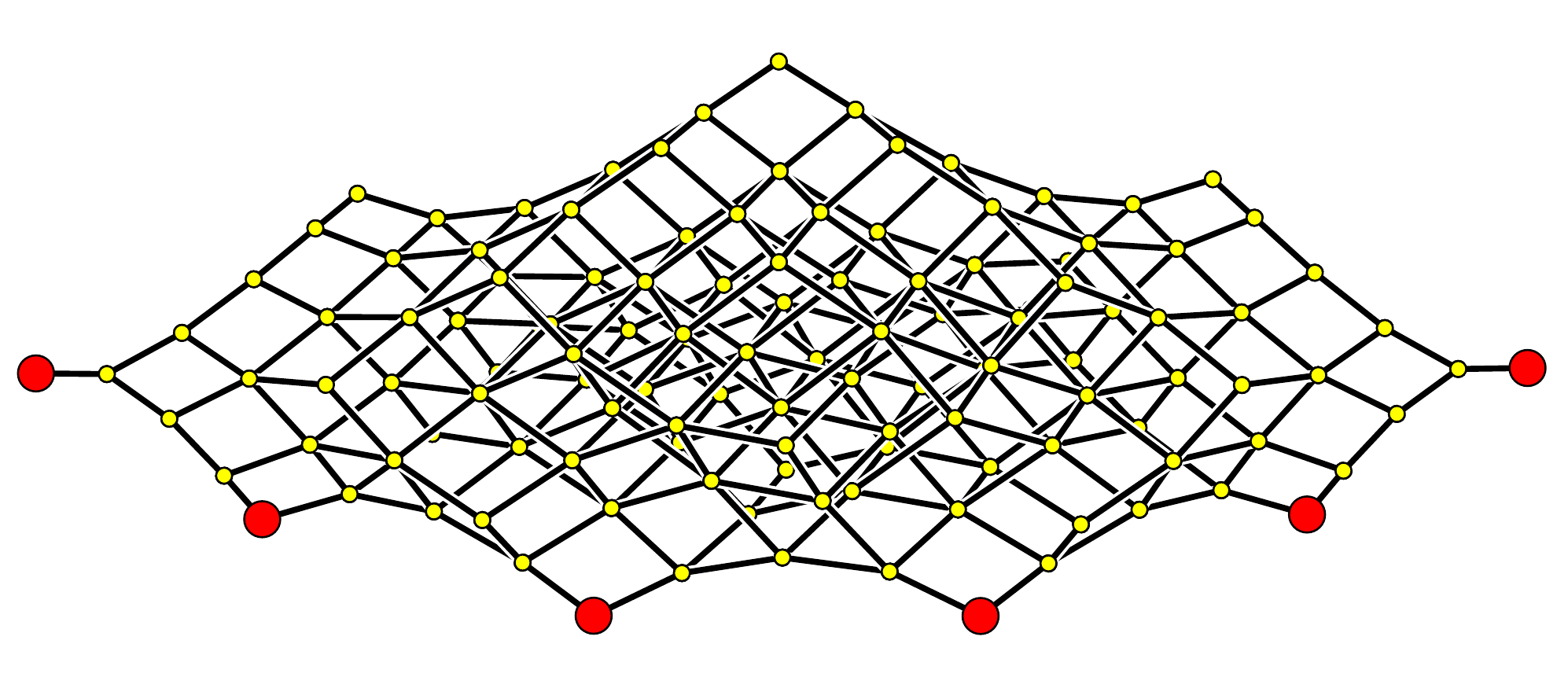}
  \caption{Tropical cyclic $4$-polytope with six tropical vertices (red).}
  \label{fig:trop-cyclic}
\end{figure}

\begin{ex} For $M(m,n)=(\mu_{ij})\in\ZZ^{m\times(n+1)}$ with $\mu_{ij}=ij$ the bounded
  subcomplex $\cT(T_{M(m,n)})$ is the \emph{tropical cyclic polytope} with $m$ vertices in
  dimension $n$ \cite{BlockYu06}.  The case $m=6$ and $n=4$ is shown in Figure
  \ref{fig:trop-cyclic}.  The tropical cyclic polytope $\cT(T_{M(4,6)})$ has $126$
  pseudo-vertices, six of which are tropical vertices.
\end{ex}

\subsection{{pd}-Graphs of Simplicial Manifolds}
\label{sec:pd_graphs}

\noindent Natural candidates for graphs to associate with a finite simplicial complex are
its $1$-skeleton (also called the \emph{primal graph}) and its dual graph.  They both
encode neighborhood information, the adjacency of vertices in the primal graph, and the
adjacency of facets in the dual graph. But each carries only incomplete information: The
primal graph does not ``see'' the facets, that is, we do not know which set of vertices of
the graph forms a face (of dimension higher than one) and which does not.  On the other
hand, we cannot tell from the dual graph if facets intersect, unless the intersection is a
ridge.  Thus it is sometimes useful to combine the two graphs in a common picture. This
allows us to find a more ``accurate'' embedding of the graphs since we may be able to make
use of the two kind of adjacency informations at the same time.  This approach will
sometimes allow for pictures of not too high-dimensional objects which carry some
geometric meaning.

We define the \emph{primal-dual graph}, or \emph{{pd}-graph} for short, as the disjoint
union of the primal and the dual graph with the following additional edges: A primal node
corresponding to a vertex~$v$ and a dual node corresponding to facet~$F$ are connected by
an \emph{artificial} edge if~$v \in F$.

\begin{figure}[htbp] \centering
  \includegraphics[width=.8\textwidth]{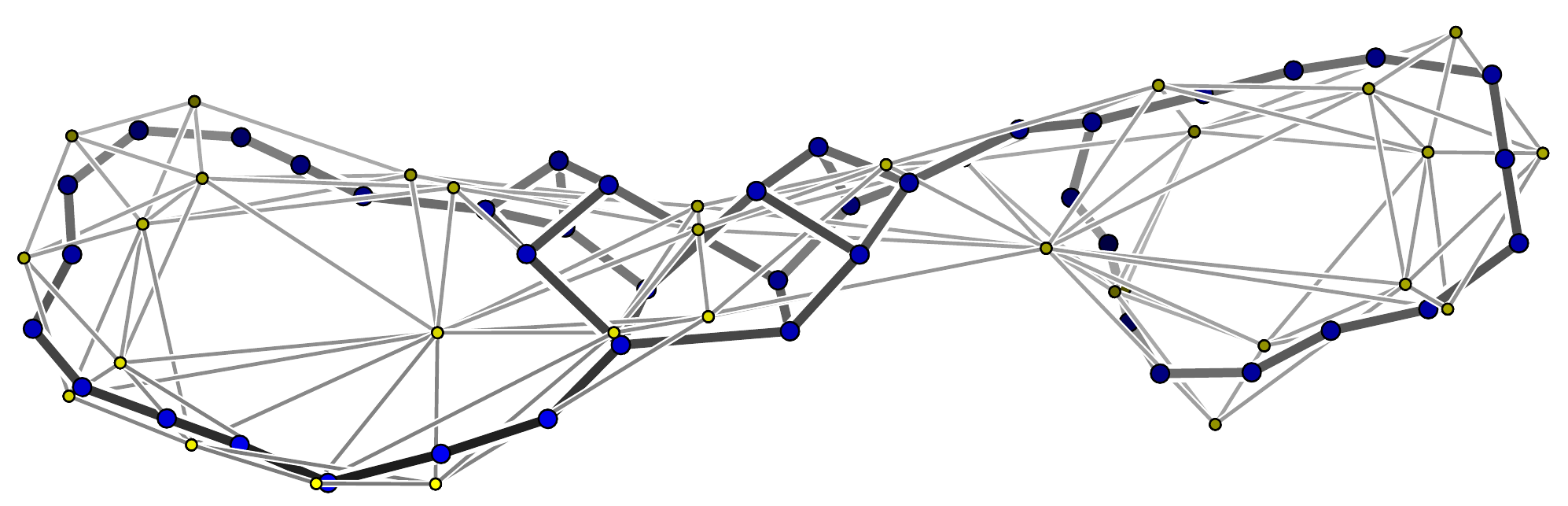}
  \caption{{pd}-graph of a triangulation of a solid surface of genus~2.  The thin edges
    belong to the primal graph, and the thick ones to the dual.}
  \label{fig:solid_g2-srf}
\end{figure}

\begin{ex} Consider an abstract simplicial complex~$K$ homeomorphic to a solid whose
  boundary surface is of genus~$2$.  A priori there are no coordinates for the vertices
  of~$K$ given, but using our pseudo-physical model (following
  Section~\ref{subsec:attr_and_repell_forces}) on the {pd}-graph of~$K$ we nevertheless
  are able to produce a decent picture, clearly showing the two holes of~$K$; see
  Figure~\ref{fig:solid_g2-srf}.  Here the artificial edges are removed to show only the
  primal and the dual graph.  If we choose the desired edge lengths of the artificial
  edges sufficiently small, then each dual node (corresponding to a facet~$F$) is pulled
  into the interior of the simplex defined by the nodes corresponding to the vertices
  of~$F$.
\end{ex}

\begin{ex} Figure~\ref{fig:min_c4} depicts the unique facet-minimal triangulation of the
  4-cube $C_4 = [0,1]^4$. The triangulation has no additional vertices and~16 facets. The
  entire $f$-vector reads $(16,57,86,60,16)$.  One way to visualize this triangulation
  of~$C_4$ is to look at the triangulation induced on the boundary of~$C_4$ in its
  Schlegel Diagram; see Section~\ref{sec:schlegel_diagrams}. Schlegel diagrams preserve
  combinatorial data as well as geometric information. But they project the polytope to
  one of its facets, hence only the boundary is visualized and all the information about
  the triangulation in the interior is lost.

\begin{figure}[htbp] \centering
  \includegraphics[width=.45\textwidth]{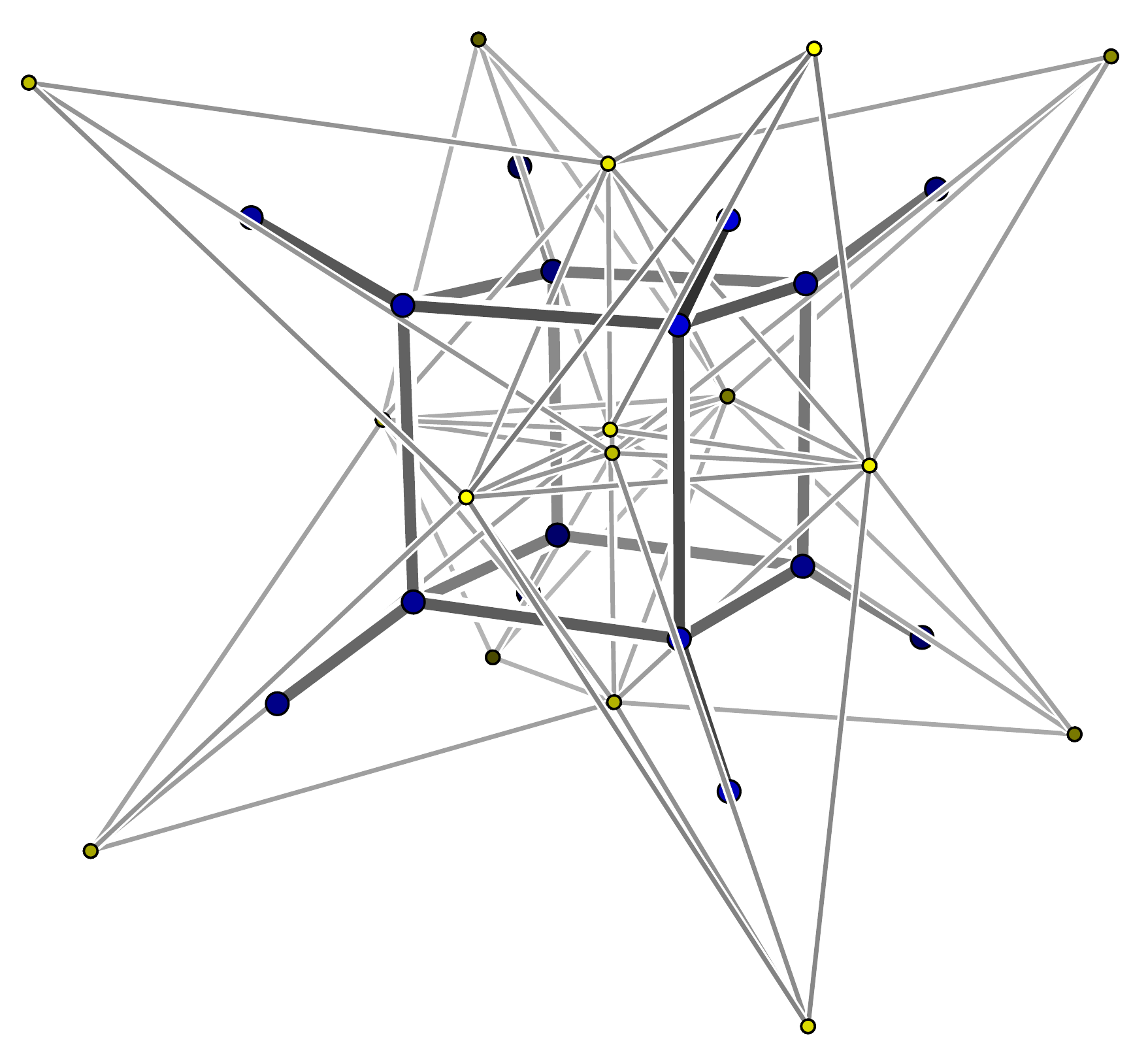}
  \caption{{pd}-graph of the unique facet-minimal triangulation of~$C_4$.}
  \label{fig:min_c4}
\end{figure}

  Alternatively, we embed the {pd}-graph of the minimal triangulation of~$C_4$ via the
  spring embedder.  Any image of a higher than 3-dimensional object has to be admired with
  care. Nevertheless one can identify the facets of the minimal triangulation of~$C_4$ and
  the way they intersect in the embedding of the {pd}-graph in Figure~\ref{fig:min_c4}.
\end{ex}

\section{Acknowledgments}

\noindent We are grateful to Ronald Wotzlaw for his beautiful regular realization of the
$4$-dimensional permutohedron which led to Figure \ref{fig:schlegel:first}.

\newpage

\bibliographystyle{amsplain} \bibliography{../bibtex/literature}

\providecommand{\bysame}{\leavevmode\hbox to3em{\hrulefill}\thinspace}
\providecommand{\MR}{\relax\ifhmode\unskip\space\fi MR }
% \MRhref is called by the amsart/book/proc definition of \MR.
\providecommand{\MRhref}[2]{%
  \href{http://www.ams.org/mathscinet-getitem?mr=#1}{#2}
}
\providecommand{\href}[2]{#2}
\begin{thebibliography}{10}

\bibitem{Balinski61}
M.~L. Balinski, \emph{On the graph structure of convex polyhedra in
  {$n$}-space}, Pacific J. Math. \textbf{11} (1961), 431--434. \MR{MR0126765
  (23 \#A4059)}

\bibitem{MR858908}
Hans-J{\"u}rgen Bandelt and Andreas Dress, \emph{Reconstructing the shape of a
  tree from observed dissimilarity data}, Adv. in Appl. Math. \textbf{7}
  (1986), no.~3, 309--343. \MR{MR858908 (87k:05060)}

\bibitem{BerEppGui-ESA-95}
Marshall~Wayne Bern, David Eppstein, Leonidas~J. Guibas, John~E. Hershberger,
  Subhash Suri, and Jan~Dithmar Wolter, \emph{{The centroid of points with
  approximate weights}}, Proc. 3rd Eur. Symp. Algorithms (ESA 1995) (Paul~G.
  Spirakis, ed.), Lecture Notes in Computer Science, no. 979, Springer-Verlag,
  September 1995, \\
  \url{http://www.ics.uci.edu/~eppstein/pubs/BerEppGui-ESA-95.ps.gz},
  pp.~460--472.

\bibitem{BlindMani87}
Roswitha Blind and Peter Mani-Levitska, \emph{Puzzles and polytope
  isomorphisms}, Aequationes Math. \textbf{34} (1987), no.~2-3, 287--297.
  \MR{MR921106 (89b:52008)}

\bibitem{BlockYu06}
Florian Block and Josephine Yu, \emph{Tropical convexity via cellular
  resolutions}, J. Algebraic Combin. \textbf{24} (2006), no.~1, 103--114.
  \MR{MR2245783 (2007f:52041)}

\bibitem{Cremona1890}
Luigi Cremona, \emph{Graphical statics}, Oxford University Press, 1890, English
  translation by {T}.~{H}.~{B}eare.

\bibitem{DevelinSturmfels04}
Mike Develin and Bernd Sturmfels, \emph{Tropical convexity}, Doc. Math.
  \textbf{9} (2004), 1--27 (electronic), correction: ibid., pp.\ 205--206.
  \MR{MR2054977 (2005i:52010)}

\bibitem{eppstein:_ukrain_easter_egg}
David Eppstein, \emph{Ukrainian easter egg}, in: ``The Geometry Junkyard'',
  computational and recreational geometry, 23 January 1997,\\
  \url{http://www.ics.uci.edu/~eppstein/junkyard/ukraine/}.

\bibitem{fruchtermann92:_graph_drawin_force_placem}
Thomas Fruchtermann and Edward Reingold, \emph{Graph drawing by force-directed
  placement}, Software Practice and Experience \textbf{21} (1992), no.~11,
  1129--1164.

\bibitem{polymake}
Ewgenij Gawrilow and Michael Joswig, \emph{{\rm\tt{polymake}}, version 2.3
  (desert)}, 1997--2007, with contributions by {T. R{\"o}rig} and {N. Witte},
  free software, \url{http://www.math.tu-berlin.de/polymake}.

\bibitem{gawrilow_joswig:POLYMAKE_I}
Ewgenij Gawrilow and Michael Joswig, \emph{{\tt polymake}: a framework for
  analyzing convex polytopes}, Polytopes--combinatorics and computation
  (Oberwolfach, 1997), DMV Sem., vol.~29, Birkh{\"a}user, Basel, 2000,
  pp.~43--73. \MR{2001f:52033}

\bibitem{JReality}
Charles Gunn, Tim Hoffmann, Markus Schmies, and Steffen Wei{\ss}mann,
  \emph{\texttt{jReality}}, \url{http://www.jreality.de}, 2007.

\bibitem{splitstree}
Daniel~H. Huson and David Bryant, \emph{Application of phylogenetic networks in
  evolutionary studies}, Mol. Biol. Evol. \textbf{23} (2006), no.~2, 254--267,
  \url{www.splitstree.org}.

\bibitem{MR0182949}
J.~R. Isbell, \emph{Six theorems about injective metric spaces}, Comment. Math.
  Helv. \textbf{39} (1964), 65--76. \MR{MR0182949 (32 \#431)}

\bibitem{eg-models_goldfarb}
Michael Joswig, \emph{Goldfarb's cube}, 2000, \url{http://www.eg-models.de},
  eg-model nr.\ 2000.09.030.

\bibitem{Joswig05}
Michael Joswig, \emph{Tropical halfspaces}, Combinatorial and computational
  geometry, Math. Sci. Res. Inst. Publ., vol.~52, Cambridge Univ. Press,
  Cambridge, 2005, pp.~409--431. \MR{MR2178330 (2006g:52012)}

\bibitem{JuengerMutzel04}
Michael J{\"u}nger and Petra Mutzel (eds.), \emph{Graph drawing software},
  Mathematics and Visualization, Springer-Verlag, Berlin, 2004. \MR{MR2159308}

\bibitem{KaibelSchwartz03}
Volker Kaibel and Alexander Schwartz, \emph{On the complexity of polytope
  isomorphism problems}, Graphs Combin. \textbf{19} (2003), no.~2, 215--230.
  \MR{MR1996205 (2004e:05125)}

\bibitem{klee-minty1972}
Victor Klee and George~J. Minty, \emph{How good is the simplex algorithm?},
  Inequalities, III (Proc. Third Sympos., Univ. California, Los Angeles,
  Calif., 1969; dedicated to the memory of Theodore S. Motzkin), Academic
  Press, New York, 1972, pp.~159--175. \MR{MR0332165 (48 \#10492)}

\bibitem{maxwell1864:on_reciprocal_figures}
J.~C. Maxwell, \emph{On reciprocal figures and diagrams of forces},
  Philosophical Magazine (1864), 250--261, Ser. 4, 27.

\bibitem{javaview}
Konrad Polthier, Klaus Hildebrandt, Eike Preuss, and Ulrich Reitebuch,
  \emph{\texttt{JavaView}, version~3.95}, \url{http://www.javaview.de}, 2007.

\bibitem{richter-gebert96:realization_spaces_of_polytopes}
J.~Richter-Gebert, \emph{Realization spaces of polytopes}, vol. 1643, Springer
  Verlag, Berlin Heidelberg, 1996.

\bibitem{roerig-witte-ziegler2007}
Thilo R{\"o}rig, Nikolaus Witte, and G{\"u}nter~M. Ziegler, \emph{Zonotopes
  with large 2d cuts}, 2007, \url{arXiv:0710.3116v2}.

\bibitem{Steinitz22}
Ernst Steinitz, \emph{{P}olyeder und {R}aumteilungen}, {E}ncyklop{\"a}die der
  mathematischen Wissenschaften, 1922, Dritter Band: Geometrie, III.1.2., Heft
  9, Kapitel 3 A B 12, pp.~1--139.

\bibitem{SteinitzRademacher76}
Ernst Steinitz and Hans Rademacher, \emph{Vorlesungen \"uber die {T}heorie der
  {P}olyeder unter {E}inschluss der {E}lemente der {T}opologie},
  Springer-Verlag, Berlin, 1976, Reprint der 1934 Auflage, Grundlehren der
  Mathematischen Wissenschaften, No. 41. \MR{MR0430958 (55 \#3962)}

\bibitem{MR2097310}
Bernd Sturmfels and Josephine Yu, \emph{Classification of six-point metrics},
  Electron. J. Combin. \textbf{11} (2004), no.~1, Research Paper 44, 16 pp.
  (electronic). \MR{MR2097310 (2005m:51016)}

\bibitem{TollisEtAl99}
Ioannis~G. Tollis, Giuseppe Di~Battista, Peter Eades, and Roberto Tamassia,
  \emph{Graph drawing}, Prentice Hall Inc., Upper Saddle River, NJ, 1999,
  Algorithms for the visualization of graphs. \MR{MR2064104 (2005i:68067)}

\bibitem{tutte63:how_to_draw_a_graph}
W.~T. Tutte, \emph{How to draw a graph}, Proc. London Math. Soc. \textbf{13}
  (1963), no.~3, 743--767.

\bibitem{OWReport}
Uli {Wagner (ed.)}, \emph{{Conference on Geometric and Topological
  Combinatorics: Problem Session}}, Oberwolfach Reports \textbf{4} (2006),
  no.~1, 265--267.

\bibitem{ziegler:LOP}
G{\"u}nter~M. Ziegler, \emph{Lectures on polytopes}, Springer, 1995.

\end{thebibliography}

\end{document}